\documentclass[a4,12pt]{amsart}
\usepackage{amssymb,amstext,amscd,amsfonts,mathdots,mathrsfs}
\numberwithin{equation}{section}
\usepackage{hyperref}
\usepackage{setspace}
\usepackage[all]{xy}
\usepackage{amsthm}
\usepackage{enumerate}
\usepackage{diagbox}
\usepackage{indentfirst}
\usepackage{color}
\usepackage{colortbl}
\usepackage{amsmath}
\usepackage{tikz}
\usetikzlibrary{shapes.geometric, arrows}
\usetikzlibrary{calc}
\usepackage{multirow,multicol}
\usepackage{listings}
\definecolor{codegreen}{rgb}{0,0.6,0}
\definecolor{codegray}{rgb}{0.5,0.5,0.5}
\definecolor{codepurple}{rgb}{0.58,0,0.82}
\definecolor{backcolour}{rgb}{0.95,0.95,0.92}
\lstdefinestyle{mystyle}{
  backgroundcolor=\color{backcolour}, commentstyle=\color{codegreen},
  keywordstyle=\color{magenta},
  numberstyle=\tiny\color{codegray},
  stringstyle=\color{codepurple},
  basicstyle=\ttfamily\footnotesize,
  breakatwhitespace=false,         
  breaklines=true,                 
  captionpos=b,                    
  keepspaces=true,                 
  numbers=left,                    
  numbersep=5pt,                  
  showspaces=false,                
  showstringspaces=false,
  showtabs=false,                  
  tabsize=2
}
\newtheorem{theorem}{Theorem}[section]
\newtheorem{theorem*}{Theorem}
\newtheorem{corollary}[theorem]{Corollary}
\newtheorem{corollary*}[theorem*]{Corollary}
\newtheorem{lemma}[theorem]{Lemma}
\newtheorem{proposition}[theorem]{Proposition}

\theoremstyle{definition}
\newtheorem{definition}[theorem]{Definition}
\newtheorem{remark}[theorem]{Remark}

\newtheorem*{question*}{Question}

\newtheorem*{conjecture*}{Conjecture}

\newtheorem*{notation*}{Notation}

\newtheorem*{claim*}{Claim}

\begin{document}
\setlength{\baselineskip}{16pt}

\title
[Representation-finite simply connected tensor products]
{Representation-finite tensor product of simply connected algebras}

\author[K. Miyamoto]{Kengo Miyamoto}
\address{Department of Computer and Information Science, Graduate School of Science and Engineering, Ibaraki University, Hitachi, Ibaraki 316-8511, Japan.}
\email{kengo.miyamoto.uz63@vc.ibaraki.ac.jp}

\author[Q. Wang]{Qi Wang}
\address{Yau Mathematical Sciences Center, Tsinghua University, Beijing 100084, China.}
\email{infinite-wang@outlook.com}


\begin{abstract}
We determine the representation-finiteness of $A\otimes B$, where both $A$ and $B$ are simply connected algebras with at least three simple modules.
\end{abstract}
\maketitle

\section{Introduction}
A fundamental theme in the representation theory of finite-dimensional algebras is classification of representation types. Let $A$ be a finite-dimensional algebra over an algebraically closed field $\mathsf{k}$. Then $A$ is said to be \emph{representation-finite} if there are only finitely many indecomposable $A$-modules. Otherwise, $A$ is said to be \emph{representation-infinite}. We say that a representation-infinite algebra $A$ is \emph{tame} if all but finitely many $d$-dimensional indecomposable $A$-modules can be organized in a one-parameter family, for each dimension $d$. A representation-infinite algebra $A$ is called \emph{wild} if there exists a faithful exact $\mathsf{k}$-linear functor from the module category of the free associative algebra $\mathsf{k}\langle x,y\rangle$ to the module category of $A$. Thanks to Drozd's Finite-Tame-Wild Trichotomy (see \cite{Dr-tame-wild}), we know that the representation type of any finite-dimensional algebra over $\mathsf{k}$ is exactly one of representation-finite, tame and wild.

The representation type of tensor product algebras has been studied in various contexts. In the 1970s, Bondarenko and Drozd \cite{BD-finite-group} considered the representation type of finite groups, while Auslander and Reiten \cite{AR-triangular-matrix} dedicated their effort to the representation type of triangular matrix rings. Moving into the 1980s, mathematicians delved into the representation type of triangular matrix algebras over different classes of algebras, for example, see \cite{S-nakayama} for Nakayama algebras, \cite{HM-selfinjective} for self-injective algebras, \cite{L-special-alg} for radical square zero algebras, etc. The most recent progress in this field can be attributed to Leszczy$\acute{\text{n}}$ski and Skowro$\acute{\text{n}}$ski, as clear in their series of papers \cite{L-rep-type, LS-tame-triangular-matrix, LS-tame-tensor-product}. These papers provide a complete description of non-wild tensor product algebras. However, it is still open to distinguish representation-finite cases and tame cases from non-wild cases.

This paper is dedicated to providing a complete classification of representation-finite tensor products between simply connected algebras. We refer to Subsection \ref{subsection-simply-connected} for the definitions of simply connected algebras.
The importance of simply connected algebras is that the indecomposable $A$-modules for any representation-finite algebra $A$ can be lifted to indecomposable $B$-modules over a simply connected algebra $B$, where $B$ is contained inside a certain Galois covering of the standard form $\widetilde{A}$ of $A$, see \cite{BG-standard form} for more details.
It is worth mentioning that Leszczy$\acute{\text{n}}$ski (\cite{L-rep-type}) has classified weakly sincere non-wild tensor product algebras, in terms of quiver and relations. In our classification, we also rely on quiver and relations; however, we do not need the restriction on weakly sincere cases.

Let $A$ and $B$ be two simply connected algebras. If $A$ has only two simple modules, then $A$ is isomorphic to the $2\times 2$ upper triangular matrix algebras, and the representation type of $A\otimes B$ is completely determined in \cite{LS-tame-triangular-matrix}. Hence, we only consider the cases in both $A$ and $B$ to have at least three simple modules, up to isomorphism. We also mention that the tensor product $A\otimes B$ is again simply connected, see Proposition \ref{prop::tensor-simply-connected}.

We denote by $\mathsf{rad}\ A$ the Jacobson radical of $A$ and by $|M|$ the number of isomorphism classes of indecomposable direct summands for an $A$-module $M$. Recall that an $A$-module $M$ is said to be \emph{uniserial} if it has a unique composition series, and $A$ is said to be \emph{Nakayama} if each indecomposable projective $A$-module and each indecomposable injective $A$-module are uniserial. Then, the main result of this paper can be presented in the following table. See Theorem \ref{theo::not-nakayama}, Theorem \ref{Theo::Nn-naka-rad-not-zero} and Theorem \ref{theo::naka-rad-zeor-not-naka}. Here, 'F' means representation-finite and 'IF' means representation-infinite.
\begin{center}
\renewcommand\arraystretch{1.5}
\scalebox{0.9}{\begin{tabular}{|c|c|c||c|c|c|c|c|c|} \hline
\multicolumn{3}{|c||}{\multirow{3}{*}{$A\otimes B$: simply connected}} &\multicolumn{3}{c|}{$B$: Nakayama} & \multicolumn{3}{c|}{\multirow{2}{*}{$B$: non-Nakayama}}\\
\cline{4-6}

\multicolumn{3}{|l||}{} &\multicolumn{2}{c|}{$\mathsf{rad}^2=0$}& \multirow{2}{*}{$\mathsf{rad}^2\neq 0$} &\multicolumn{3}{l|}{}  \\
\cline{4-5} \cline{7-9}

\multicolumn{3}{|l||}{} &$|B|=3$&$|B|\geq 4$&&$|B|=3$&$|B|=4$&$|B|\geq 5$\\ \hline\hline

\multirow{3}*{$A$: Nakayama}&\multirow{2}*{$\mathsf{rad}^2=0 $}&$|A|=3$&\multicolumn{2}{c|}{\multirow{2}{*}{F}}&\multirow{2}*{F$\&$IF}&\multirow{2}*{F}&\multirow{2}*{F$\&$IF}&F$\&$IF\\
\cline{3-3} \cline{9-9}

 &  & $|A|\geq 4$ & \multicolumn{2}{c|}{ } & & & &IF\\
\cline{2-9}

& \multicolumn{2}{c||}{$\mathsf{rad}^2\neq 0$} &\multicolumn{2}{c|}{F$\&$IF}&IF& \multicolumn{3}{c|}{IF}\\ \hline

\multicolumn{2}{|c|}{\multirow{3}{*}{$A$: non-Nakayama}}  &$|A|=3$& \multicolumn{2}{c|}{F} &\multirow{3}*{IF} & \multicolumn{3}{c|}{ }  \\
\cline{3-5}

\multicolumn{2}{|c|}{}&$|A|=4$&\multicolumn{2}{c|}{F$\&$IF}& &\multicolumn{3}{c|}{IF}\\
\cline{3-5}

\multicolumn{2}{|c|}{}&$|A|\geq 5$&F$\&$IF&IF& &\multicolumn{3}{c|}{ }  \\ \hline
\end{tabular}}
\end{center}

An immediate corollary of our results is that we have obtained a complete classification of $\tau$-tilting finite simply connected tensor products. This follows from the work of the second author as demonstrated in \cite{W-simply}, where it is shown that a simply connected algebra is $\tau$-tilting finite if and only if it is representation-finite. Here, $\tau$-tilting finiteness is a modern analog of representation-finiteness, which is inspired by the $\tau$-tilting theory initially proposed by Adachi, Iyama and Reiten \cite{AIR}.

This paper is organized as follows. In Section 2, we review some basic materials and fix our notations, including simply connected algebras, tame concealed algebras, Tits form of algebras, etc. In Section 3, we consider the case that one of $A$ and $B$ is a representation-finite hereditary algebra, and the results here will reduce the problem to the case that one of $A$ and $B$ is a Nakayama algebra with radical square zero, which will be handled in Section 4. Finally, we give a complete list of minimal representation-infinite simply connected tensor product algebras in Section 5.

\section{Preliminaries}
Let $A=\mathsf{k}Q/\mathcal{I}$ be a bound quiver algebra with a finite connected quiver $Q=(Q_0, Q_1)$ and an admissible ideal $\mathcal{I}$ over an algebraically closed field $\mathsf{k}$.
A relation $\rho=\sum_{i=1}^{m}\lambda_i\omega_i$ in $\mathcal{I}$ is a $\mathsf{k}$-linear combination of paths $\omega_i$ of length at least two having the same source and target, where the $\lambda_i$ are scalars and not all zero. If $m=2$, then $\rho$ is called a \emph{commutativity relation}. Throughout, we use a dotted line to indicate the obvious commutativity relation in the algebra. For example, the dotted lines in the following quiver indicate $ab=ba$:
$$
\vcenter{\xymatrix@C=1.5cm@R=0.7cm{
\circ \ar[r]^a\ar[d]^b\ar@{.}[dr]&\circ \ar[d]^b &\circ\ar[l]_a\ar[d]^b\ar@{.}[dl]&\circ\ar[l]_a \ar[r]^a\ar[d]^b\ar@{.}[dr]\ar@{.}[dl]&\circ\ar[d]^b\\
\circ \ar[r]^a&\circ  &\circ\ar[l]_a&\circ\ar[l]_a \ar[r]^a&\circ
}}.
$$
A relation $\rho\in\mathcal{I}$ is said to be \emph{minimal} if $m\geq 2$, and $\sum_{i\in J}\lambda_iw_i\not\in \mathcal{I}$ for any non-empty proper subset $J\subset\{1,2,\ldots,m\}$. Moreover, we sometimes use $\rho=0$ to indicate $\rho\in\mathcal{I}$.
We refer to \cite{ASS} for more background on the representation theory of quivers.

A full subquiver $Q'$ of $Q$ is said to be \emph{convex} if, for any path $a_0\to a_1\to\cdots\to a_m$, $a_0, a_m\in Q'$ implies $a_i\in Q'$ for all $1\leq i\leq m$.
By a \emph{convex subalgebra} of $A=\mathsf{k}Q/\mathcal{I}$, we mean a bound quiver algebra $\mathsf{k}Q'/\mathcal{I}'$ with a convex subquiver $Q'$ of $Q$ and $\mathcal{I}'=\mathcal{I}\cap\mathsf{k}Q'$.

We denote by $\mathsf{rad}\ A$ the Jacobson radical of $A$ and by $A^{\mathsf{op}}$ the opposite algebra of $A$. Then, an algebra $A$ is said to be \emph{radical square zero} if $\mathsf{rad}^2\ A=0$. Recall that an algebra $A=\mathsf{k}Q/\mathcal{I}$ is called \emph{triangular} if $Q$ does not have loops and oriented cycles.

\subsection{Tensor product of algebras}
Let $A=\mathsf{k}Q_A/\mathcal{I}_A$ and $B=\mathsf{k}Q_B/\mathcal{I}_B$ be two bound quiver algebras. The tensor product $A\otimes B$ is a $\mathsf{k}$-algebra defined by the multiplication $(a_1\otimes b_1)(a_2\otimes b_2)=a_1a_2\otimes b_1 b_2$, for $a_1, a_2\in A$ and $b_1, b_2\in B$. The quiver presentation of $A\otimes B$ is given as follows.
\begin{itemize}
\item Let $Q_A\otimes Q_B$ be the quiver consisting of the vertex set $(Q_A\otimes Q_B)_0=(Q_A)_0\times (Q_B)_0$
and the arrow set $(Q_A\otimes Q_B)_1=((Q_A)_1\times (Q_B)_0)\cup ((Q_A)_0\times (Q_B)_1)$, together with the source map $s(-)$ and the target map $t(-)$ defined as
\begin{align*}
&s(\alpha\times v)=s_A(\alpha)\times v, \quad  s(u\times \beta)=u\times s_B(\beta),  \\
 &t(\alpha\times v)=t_A(\alpha)\times v, \quad t(u\times \beta)=u\times t_B(\beta),
\end{align*}
for $(\alpha,v)\in (Q_A)_1\times (Q_B)_0$ and $(u,\beta)\in (Q_A)_0\times (Q_B)_1$.

\item Let $\mathcal{I}_A\diamond \mathcal{I}_B$ be the two-sided ideal of $\mathsf{k}(Q_{A}\otimes Q_B)$ generated by the elements in $((Q_A)_0 \times \mathcal{I}_B)\cup (\mathcal{I}_A\times (Q_B)_0)$ and the elements being of the form $(a,\beta_{cd})(\alpha_{ab},d)-(\alpha_{ab},c)(b,\beta_{cd})$, where $\alpha_{ab}$ and $\beta_{cd}$ run through all arrows $\alpha_{ab}:a\rightarrow b$ in $(Q_A)_1$ and $\beta_{cd}:c\rightarrow d$ in $(Q_B)_1$.
\end{itemize}
We have $A\otimes B\simeq \mathsf{k}(Q_A\otimes Q_B)/ (\mathcal{I}_A\diamond\mathcal{I}_B)$, see \cite[Lemma 1.3]{L-rep-type}.

Some basic properties of tensor products are listed as follows.
\begin{proposition}[{\cite[Lemma 1.2]{L-rep-type}}]
Let $A$, $B$ and $C$ be bound quiver algebras. Then,
\begin{enumerate}
\item[\rm{(1)}] $A\otimes B \simeq B\otimes A$.
\item[\rm{(2)}] $(A\otimes B)^{\mathsf{op}} \simeq A^{\mathsf{op}} \otimes B^{\mathsf{op}} $.
\item[\rm{(3)}] $A\otimes (B \otimes C)\simeq (A\otimes B)\otimes C$.
\item[\rm{(4)}] $A\otimes B$ is connected if and only if both $A$ and $B$ are connected.
\end{enumerate}
\end{proposition}

We give a new property of tensor products here. Let $\{e_1,e_2, \ldots, e_n\}$ be a complete list of pairwise orthogonal primitive idempotents of $A$. Then, $A$ is said to be \emph{Schurian} if $\mathsf{dim}_{\mathsf{k}}(e_iAe_j) \leq 1$, for all $1\le i,j\le n$.
\begin{proposition}
Let $A$ and $B$ be two triangular algebras. Then, $A\otimes B$ is Schurian if and only if both $A$ and $B$ are Schurian.
\end{proposition}
\begin{proof}
Let $\{e_1^A,\ldots ,e_n^A\}$ and $\{e_1^B,\ldots, e_m^B\}$ be complete sets of pairwise orthogonal primitive idempotents of $A$ and $B$, respectively. Then, $\{e_i^A\otimes e_j^B\ | \ 1\leq i\leq n,\ 1\leq j\leq m\}$ is a complete set of pairwise orthogonal primitive idempotents of $A\otimes B$. We have
$$
\mathsf{dim}_{\mathsf{k}}(e_i^A\otimes e_j^B)(A\otimes B)(e_k^A\otimes e_l^B) = \mathsf{dim}_{\mathsf{k}}(e_i^AAe_k^A) \times \mathsf{dim}_{\mathsf{k}}(e_j^BBe_l^B). 
$$
It implies that $A$ and $B$ are Schurian if and only if $A\otimes B$ is Schurian.
\end{proof}

\subsection{Simply connected algebra}\label{subsection-simply-connected}
We recall the basics of simply connected algebras. 
See \cite{Assem-simply-connected, AS, MP} for more materials.

Let $A=\mathsf{k}Q/\mathcal{I}$ be a triangular algebra with a connected quiver $Q=(Q_0,Q_1)$ and an admissible ideal $\mathcal{I}$. For an arrow $\alpha\in Q_1$, we write $\alpha^{-1}$ as the formal inverse of $\alpha$. Let $a$ and $b$ be two vertices of $Q$. A walk from $a$ to $b$ is defined as a formal composition $\alpha_1^{\varepsilon_1}\alpha_2^{\varepsilon_2}\cdots\alpha_m^{\varepsilon_m}$, where $\alpha_i\in Q_1$ and $\varepsilon_i\in\{\pm 1\}$ for $i=1,2,\ldots ,m$.
For each vertex $a\in Q_0$, we regard the trivial path $e_a$ as the stationary walk at $a$. If $w$ is a walk from $a$ to $b$ and $w'$ is a walk from $b$ to $c$, the multiplication $ww'$ is given by concatenation of $w$ and $w'$.
We denote by $Q^{\ast}$ the set of all walks of $Q$.
Then, the homotopy relation $\sim _{\mathcal{I}}$ in $Q^{\ast}$ is defined to be the smallest equivalence relation in $Q^{\ast}$ satisfying the following conditions:
\begin{itemize}
\item $\alpha\alpha^{-1}\sim _{\mathcal{I}} e_a$ and $\alpha^{-1}\alpha\sim _{\mathcal{I}} e_b$, for each arrow $a\xrightarrow{\alpha} b$.
\item For each minimal relation $\sum_{i=1}^m\lambda_iw_i$ in $\mathcal{I}$, we have $w_i\sim _{\mathcal{I}} w_j$ for all $1\leq i,j\leq m$.
\item If $u,v,w$ and $w'$ are walks such that $u\sim _{\mathcal{I}} v$ and $w\sim _{\mathcal{I}} w'$, then we have $wuw'\sim _{\mathcal{I}} wvw'$, whenever the multiplications are defined.
\end{itemize}
We denoted by $[w]$ the equivalence class of a walk $w$ in $Q^{\ast}$. We define $\pi_1(Q,\mathcal{I},a)$ as the set of equivalence classes of all walks from $a$ to $a$. It is easily seen that $\pi_1(Q,\mathcal{I},a)$ becomes a group via the multiplication $[w]\cdot[w']=[ww']$ induced by the multiplication on $Q^{\ast}$. It is also not difficult to find that $\pi_1(Q,\mathcal{I},a)$ does not depend on the choice of $a\in Q_0$. We call the group $\pi_1(Q,\mathcal{I}):=\pi_1(Q,\mathcal{I},a)$ the \emph{fundamental group} of $A$.

\begin{definition}[{\cite[Definition 1.2]{AS}}]\label{def::simply-connected}
A triangular $\mathsf{k}$-algebra $A$ is said to be simply connected if, for any presentation $A\simeq \mathsf{k}Q/\mathcal{I}$ as a bound quiver algebra, the fundamental group $\pi_1(Q,\mathcal{I})$ is trivial.
\end{definition}

In this paper, the following statement is crucial.
\begin{proposition}[{\cite[Lemma 1.7]{L-rep-type}}]\label{prop::tensor-simply-connected}
Let $A$ and $B$ be two triangular algebras. Then, $A\otimes B$ is simply connected if and only if both $A$ and $B$ are simply connected.
\end{proposition}

We give a characterization of simply connected Nakayama algebras as follows.
\begin{proposition}[{\cite[V.3.2]{ASS}}]
A simply connected algebra $A$ is Nakayama if and only if its Gabriel quiver is
$$
\overset{\rightarrow}{A_n}: \xymatrix@C=0.7cm@R=0.5cm{
1\ar[r]&2\ar[r]&3\ar[r]&4\ar[r]&\cdots\ar[r]& n-1\ar[r]& n},
$$
or equivalently, if and only if its Gabriel quiver does not contain one of $\xymatrix@C=0.5cm{\circ&\circ\ar[l]\ar[r]&\circ}$ and $\xymatrix@C=0.5cm{\circ\ar[r]&\circ &\circ\ar[l]}$ as a subquiver.
\end{proposition}

\begin{definition}
We define $N(n):=\mathsf{k}\overset{\rightarrow}{A_n}/(\mathsf{rad}^2\ \mathsf{k}\overset{\rightarrow}{A_n})$.
\end{definition}

\subsection{Tame concealed algebra}
Let $\mathsf{mod}\ A$ be the category of finitely generated right $A$-modules. We denote by $\Gamma (\mathsf{mod}\ A)$ the Auslander-Reiten quiver of $A$. A connected component $\mathcal{C}$ of $\Gamma (\mathsf{mod}\ A)$ is called \emph{preprojective} if there is no oriented cycle in $\mathcal{C}$, and any module in $\mathcal{C}$ is of the form $\tau^{-n}(P)$ for some $n\in\mathbb{N}$ and an indecomposable projective $A$-module $P$, where $\tau$ is the Auslander-Reiten translation.

Recall that an $A$-module $T$ is called a \emph{tilting module} (see \cite{HR-tilted}) if $|T|=|A|$, $\mathsf{Ext}_A^1(T,T)=0$ and the projective dimension of $T$ is at most one. Let $A:=\mathsf{k}Q$ be a hereditary algebra and $T$ a tilting $A$-module. We call the endomorphism algebra $B:=\mathsf{End}_A\ T$ a \emph{tilted algebra} of type $Q$. If moreover, any indecomposable direct summand of $T$ is contained in a preprojective component $\mathcal{C}$ of $\Gamma (\mathsf{mod}\ A)$, we call $B$ a \emph{concealed algebra} of type $Q$.
\begin{definition}
A tame concealed algebra is one of the concealed algebras of Euclidean type $\widetilde{\mathbb{A}}_n$, $\widetilde{\mathbb{D}}_n (n\ge 4)$, $\widetilde{\mathbb{E}}_6$, $\widetilde{\mathbb{E}}_7$ and $\widetilde{\mathbb{E}}_8$.
\end{definition}

Recall that an algebra $A$ is called a \emph{minimal algebra of infinite representation type} if $A$ is representation-infinite, but $A/AeA$ is representation-finite for any non-zero idempotent $e$ of $A$. Then, we have the following statement.

\begin{proposition}[{\cite[Theorem 2]{HV}}]
An algebra $A$ is a minimal algebra of infinite representation type with a preprojective component if and only if it is a tame concealed algebra or a generalized Kronecker algebra.
\end{proposition}

It is shown in \cite[(1.2)]{AS} that tame concealed algebras are simply connected. We also mention that tame concealed algebras have been completely classified by quiver and relations in \cite{HV} (see also \cite{B-critical}). By looking in detail at the shape of the quivers in \cite{HV}, one may find the following interesting fact.
\begin{remark}\label{rem::shape-4*8}
Let $A$ be a concealed algebra of type $\widetilde{\mathbb{E}}_6$, $\widetilde{\mathbb{E}}_7$ and $\widetilde{\mathbb{E}}_8$. Then, the vertices of the quiver of $A$ can be lied in at most 4 rows and 8 columns.
\end{remark}

\subsection{Tits form}
Let $A= \mathsf{k}Q/\mathcal{I}$ be a triangular algebra and ${\bf{N}}:=\{0,1,2,\dots\}$. We recall from \cite[Section 2]{B-tits-form} that the \emph{Tits form} $\mathbf{q}_A: \mathbb{Z}^{Q_0}\rightarrow\mathbb{Z}$ of $A$ is the integral quadratic form defined by
\begin{center}
$\mathbf{q}_A(v):=\sum\limits_{i\in Q_0}v_i^2-\sum\limits_{(i\rightarrow j) \in Q_1}v_iv_j+\sum\limits_{i,j\in Q_0}r(i,j)v_iv_j$,
\end{center}
where $v:=(v_i)\in \mathbb{Z}^{Q_0}$ and $r(i,j)=\left | R\cap \mathcal{I} \right |$ is the number of relations with source $i$ and target $j$, for a minimal set $R$ of relations in $\mathcal{I}$. Then, the Tits form $\mathbf{q}_A$ is called \emph{weakly positive} if $\mathbf{q}_A(v)>0$ for any $v\neq 0$ in ${\bf{N}}^{Q_0}$. We mention the following efficient method to determine the representation-finiteness of simply connected algebras.
\begin{proposition}[{\cite[XX, Theorem 2.9]{SS2}}]\label{prop::Tits-form}
Let $A$ be a simply connected algebra. Then, $A$ is representation-finite if and only if the Tits form $\mathbf{q}_A$ is weakly positive, or equivalently, if and only if $A$ does not contain a convex subalgebra which is a concealed algebra of type $\widetilde{\mathbb{D}}_n$ $(n\geq 4)$, $\widetilde{\mathbb{E}}_6$, $\widetilde{\mathbb{E}}_7$ or $\widetilde{\mathbb{E}}_8$.
\end{proposition}

\subsection{Reduction theorems}
We review some general reduction methods on the representation-finiteness of algebras, especially, tensor product algebras. The following well-known statements are due to the existence of some full faithful functors or algebra isomorphisms.
\begin{proposition}
If $A$ is a representation-finite algebra, then
\begin{enumerate}
\item[\rm{(1)}] the quotient $A/\mathcal{J}$ is representation-finite for any two-sided ideal $\mathcal{J}$ of $A$.
\item[\rm{(2)}] the idempotent truncation $eAe$ is representation-finite for any idempotent $e$ of $A$.
\item[\rm{(3)}] the opposite algebra $A^{\mathsf{op}}$ is representation-finite.
\end{enumerate}
\end{proposition}

\begin{corollary}
If $A\otimes B$ is representation-finite, then $A$ and $B$ are representation-finite.
\end{corollary}

An algebra $A$ is said to be \emph{local} if it has only one simple module, up to isomorphism. Otherwise, $A$ is called a \emph{non-local} algebra.
\begin{proposition}[{\cite[Proposition 4.1]{Aihara-Honma}}] \label{prop::three-tensor}
If $A, B$ and $C$ are non-local algebras, then $A\otimes B\otimes C$ is representation-infinite.
\end{proposition}
The proof of Proposition \ref{prop::three-tensor} is quite simple. One may observe a tame concealed algebra of type $\widetilde{\mathbb{A}}_5$ as a quotient of $A\otimes B\otimes C$, i.e., the black points in the following quiver
$$
\vcenter{\xymatrix@C=0.3cm@R=0.2cm{& \circ \ar[rr]\ar'[d][dd] \ar[dl]& & \bullet \ar[dd]\ar[dl]\\ \bullet \ar[rr]\ar[dd]& & \bullet \ar[dd]\\& \bullet \ar[dl]\ar'[r][rr]& & \bullet \ar[dl]\\ \bullet\ar[rr]& & \circ}}.
$$

Lastly, we give a new observation as follows. This may be regarded as a partial reason why we study tensor products between simply connected algebras as a first step.
\begin{proposition}
Let $A=\mathsf{k}Q/\mathcal{I}$ be a bound quiver algebra. If $Q$ contains an oriented cycle that is not a loop, then $A\otimes A$ is representation-infinite.
\end{proposition}
\begin{proof}
It suffices to show that $\Lambda:=(A\otimes A)/\mathsf{rad}^2(A\otimes A)$ is representation-infinite. Suppose $Q$ is a cycle of the form
\begin{center}
$\begin{xy}
          (0,-5)*[o]+{1}="1",(15,-5)*[o]+{2}="2",
          (30,-5)*[o]+{3}="3",(45,-5)*[o]+{\cdots}="4",
          (65,-5)*[o]+{n-3}="5",(85,-5)*[o]+{n-2}="6",
          (105,-5)*[o]+{n-1.}="7",(50,5)*[o]+{n}="8"
          \ar "1";"2"
          \ar "2";"3"
          \ar "3";"4"
          \ar "4";"5"
          \ar "5";"6"
          \ar "6";"7"
          \ar "7";"8"
          \ar "8";"1"
 \end{xy}$
\end{center}
By the construction of $\Lambda$, the quiver of $\Lambda$ contains two kinds of cycles for any $1\leq k\leq n$:
\begin{center}
$\begin{xy}
          (0,-5)*[o]+{(k,1)}="1",(20,-5)*[o]+{(k,2)}="2",
          (40,-5)*[o]+{\cdots}="3",(60,-5)*[o]+{(k,n-2)}="4",
          (85,-5)*[o]+{(k,n-1),}="5",(42,5)*[o]+{(k,n)}="6"
          \ar "1";"2"
          \ar "2";"3"
          \ar "3";"4"
          \ar "4";"5"
          \ar "5";"6"
          \ar "6";"1"
 \end{xy}$
\end{center}
and
\begin{center}
$\begin{xy}
          (0,-5)*[o]+{(1,k)}="1",(20,-5)*[o]+{(2,k)}="2",
          (40,-5)*[o]+{\cdots}="3",(60,-5)*[o]+{(n-2,k)}="4",
          (85,-5)*[o]+{(n-1,k).}="5",(42,5)*[o]+{(n,k)}="6"
          \ar "1";"2"
          \ar "2";"3"
          \ar "3";"4"
          \ar "4";"5"
          \ar "5";"6"
          \ar "6";"1"
 \end{xy}$
\end{center}
Then, $\Lambda$ admits the following quiver as a subquiver:
\begin{center}
$\begin{xy}
           (0,0)*[o]+{(1,1,0)}="1",(20,10)*[o]+{(2,1,1)}="2",
          (45,10)*[o]+{(2,n,0)}="3",(70,10)*[o]+{(3,n,1)}="4",
          (90,10)*[o]+{\cdots}="d1",(90,-10)*[o]+{\cdots}="d2",
          (70,-10)*[o]+{(n,3,1)}="6",
          (45,-10)*[o]+{(n,2,0)}="7",(20,-10)*[o]+{(1,2,1)}="8",
          (107,10)*[o]+{(s+1,s+2,\epsilon)}="k-1",(107,-10)*[o]+{(s+2,s+1,\epsilon)}="k-2",
          (135,0)*[o]+{(s+1,s+1,\epsilon)}="k",
          \ar "1";"2"
          \ar "1";"8"
          \ar "3";"2"
          \ar "3";"4"
          \ar "d1";"4"
          \ar "d2";"6"
          \ar "7";"6"
          \ar "7";"8"
          \ar @{-}"k-1";"k"^{\epsilon}
          \ar @{-}"k-2";"k"^{\epsilon}
 \end{xy}$
\end{center}
where
\begin{center}
$n=\left\{\begin{array}{ll}
 2s-1 & \text{if\ } n \text{\ is odd,} \\
 2s   & \text{if\ } n \text{\ is even,} \end{array}\right.$
$\epsilon =\left\{
\begin{array}{ll}
 1 & \text{if\ } n \text{\ is odd,} \\
 0   & \text{if\ } n \text{\ is even,} \end{array}\right.$
$\begin{xy}
(0,0)*[o]+{}="1",(10,0)*[o]+{}="2",
\ar @{-}"1";"2"^{\epsilon}
\end{xy}$
$=\left\{\begin{array}{ll}
\xrightarrow{\ \epsilon\ } & \text{if\ } n \text{\ is odd,} \\
\xleftarrow{\ \epsilon\ }  & \text{if\ } n \text{\ is even.}
\end{array}\right.$
\end{center}
The above subquiver indicates a tame concealed algebra of type $\widetilde{\mathbb{A}}_n$. Hence, $\Lambda$ is representation-infinite.
\end{proof}

\section{Hereditary algebras}
In this section, we consider some elementary cases (i.e., $A\otimes B$ with $A$ being a non-local representation-finite hereditary algebra) to deduce that most tensor products are representation-infinite. According to Gabriel's Theorem, a hereditary algebra $A$ is representation-finite if and only if the underlying graph of its quiver $Q_A$ is one of Dynkin diagrams.
We then consider the path algebras of type $A, D, E$ as follows.

\subsection{Path algebra of type $A$}
Let $\mathbb{A}_n$ ($n\geq 2$) be the Dynkin diagram of type $A$:
$$
\xymatrix@C=1cm@R=0.5cm{
1\ar@{-}[r]&2\ar@{-}[r]&3\ar@{-}[r]&4\ar@{-}[r]&\cdots\ar@{-}[r]& n-1\ar@{-}[r]& n}.
$$
The orientation of $\mathbb{A}_n$ is defined as $\varepsilon:=(\varepsilon_1,\varepsilon_2,\ldots, \varepsilon_{n-1})$ with
$$
\left\{\begin{array}{ll}
\varepsilon _i=+ & \text{if } i \longrightarrow i+1 , \\
\varepsilon _i=- & \text{if } i+1 \longrightarrow i .
\end{array}\right.
$$
We denote by $\mathbf{A}_n^\varepsilon$ the path algebra of type $A$ associated with the orientation $\varepsilon$.
Note that, if $\varepsilon_i=+$ (or $-$) for all $1\le i\le n-1$, the tensor product $\mathbf{A}_n^{\varepsilon}\otimes B$ is isomorphic to the lower triangular matrix algebra $T_n(B)$, that is,
$$
T_n(B)=T_n(K)\otimes B:=\left (\begin{smallmatrix}
B&0&\cdots &0 \\
B&B&\cdots&0\\
\vdots  &\vdots   &\ddots   &\vdots  \\
B&B&\cdots   &B
\end{smallmatrix}\right ).
$$

\begin{proposition}\label{prop::A2}
Let $B$ be a non-local representation-finite simply connected algebra. Then the following statements are equivalent.
\begin{enumerate}
\item[\rm{(1)}] $\mathbf{A}_2^{(+)}\otimes B$ is representation-finite.
\item[\rm{(2)}] $B$ and $B^\mathsf{op}$ do not contain one of the algebras in the family (IT) in \cite{LS-tame-triangular-matrix} as a quotient algebra.
\end{enumerate}
\end{proposition}
\begin{proof}
Let $C$ be a representation-finite standard algebra in the sense of \cite[(3.1)]{BG-standard form}. Then, it is shown in \cite[Theorem 4]{LS-tame-triangular-matrix} that $\mathbf{A}_2^{(+)}\otimes C$ is representation-finite if and only if its Galois covering $\widetilde{C}$ (and $\widetilde{C}^\mathsf{op}$) do not contain one of the algebras in the family (IT) as a quotient algebra. On the other hand, it is known from \cite[(6.1)]{BG-covering} that any representation-finite simply connected algebra $B$ is standard and $B\simeq \widetilde{B}$.
\end{proof}

\begin{corollary}\label{cor::A2-An}
Let $\varepsilon$ be an orientation of $\mathbb{A}_n$. Then, the tensor product $\mathbf{A}_2^{(+)}\otimes \mathbf{A}_n^{\varepsilon}$ is representation-finite if and only if $n\leq 4$.
\end{corollary}
\begin{proof}
This is obvious since $\mathbb{A}_4$ is not in (IT) but $\mathbb{A}_5$ is in (IT), see \cite[Page 277]{LS-tame-triangular-matrix}.
\end{proof}

\begin{lemma}\label{lem::A3-A3}
The tensor product $\mathbf{A}_3^\varepsilon\otimes \mathbf{A}_3^\omega$ is representation-infinite for any $\varepsilon$ and $\omega$.
\end{lemma}
\begin{proof}
Up to isomorphism, it suffices to consider the following 4 cases. We shall prove that the tensor product $\mathbf{A}_3^\varepsilon\otimes \mathbf{A}_3^\omega$ in each case admits a tame concealed algebra as a quotient, which is indicated by the black points in the corresponding quivers.
\begin{itemize}
\item Set $\varepsilon=(++)$ and $\omega=(++)$. Then, $\mathbf{A}_3^{(++)}\otimes \mathbf{A}_3^{(++)}$ is presented as
$$
\vcenter{\xymatrix@C=1.2cm@R=0.5cm{\circ \ar[r] \ar[d]\ar@{.}[dr]&\bullet\ar@{.}[dr]\ar[r] \ar[d] &\bullet\ar[d]\\ \bullet \ar@{.}[dr]\ar[r] \ar[d]&\bullet\ar@{.}[dr]\ar[r] \ar[d]&\bullet \ar[d]\\ \bullet\ar[r]&\bullet  \ar[r] &\circ}},
$$
and it admits a tame concealed algebra of type $\widetilde{\mathbb{D}}_6$ as a quotient algebra.

\item Set $\varepsilon=(++)$ and $\omega=(+-)$. Then, $\mathbf{A}_3^{(++)}\otimes \mathbf{A}_3^{(+-)}$ is presented as
$$
\vcenter{\xymatrix@C=1.2cm@R=0.5cm{\circ \ar[r] \ar[d]\ar@{.}[dr]&\bullet\ar@{.}[dr]\ar[r] \ar[d] &\bullet\ar[d]\\ \bullet \ar[r] &\bullet\ar[r] &\bullet \\ \circ\ar[r]\ar[u]\ar@{.}[ur]&\bullet\ar@{.}[ur]\ar[r] \ar[u]&\bullet  \ar[u]}},
$$
and it has a tame concealed algebra of type $\widetilde{\mathbb{E}}_6$ (see No.85 in \cite{B-critical}) as a quotient.

\item Set $\varepsilon=(+-)$ and $\omega=(+-)$. Then, $\mathbf{A}_3^{(+-)}\otimes \mathbf{A}_3^{(+-)}$ is presented as
$$
\vcenter{\xymatrix@C=1.2cm@R=0.5cm{\circ \ar[r] \ar[d]\ar@{.}[dr]&\bullet\ar[d] &\circ \ar@{.}[dl] \ar[l]\ar[d]\\ \bullet \ar[r] &\bullet &\bullet\ar[l] \\ \circ\ar[r]\ar[u]\ar@{.}[ur]&\bullet \ar[u]&\circ   \ar@{.}[ul] \ar[l]\ar[u]}},
$$
and it admits a tame concealed algebra of type $\widetilde{\mathbb{D}}_4$ as a quotient algebra.

\item Set $\varepsilon=(+-)$ and $\omega=(-+)$. Then, $\mathbf{A}_3^{(+-)}\otimes \mathbf{A}_3^{(-+)}$ is presented as
$$
\vcenter{\xymatrix@C=1.2cm@R=0.5cm{\bullet \ar[r]&\circ &\bullet \ar[l]\\ \bullet \ar@{.}[ur]  \ar@{.}[dr]\ar[r]\ar[u]\ar[d] &\bullet \ar[u]\ar[d]&\bullet  \ar@{.}[dl]\ar@{.}[ul]\ar[u]\ar[d]\ar[l] \\ \bullet\ar[r]&\circ &\bullet\ar[l]}},
$$
and it admits a tame concealed algebra of type $\widetilde{\mathbb{D}}_6$ as a quotient algebra.
\end{itemize}
This completes the proof.
\end{proof}

\begin{lemma}\label{lem::A3}
Let $B$ be a non-local representation-finite simply connected algebra. Then, $\mathbf{A}_3^{\varepsilon}\otimes B$ is representation-finite if and only if $B$ is isomorphic to $\mathbf{A}_2^{(+)}$ or $N(n)$ for $n\ge 3$.
\end{lemma}
\begin{proof}
If $B$ contains $\mathbf{A}_3^{\omega}$ for some $\omega$ as a quotient, then there exists a surjection $\mathbf{A}_3^{\varepsilon}\otimes B\to \mathbf{A}_3^{\varepsilon}\otimes \mathbf{A}_3^{\omega}$.
It follows from Lemma \ref{lem::A3-A3} that $\mathbf{A}_3^{\varepsilon}\otimes B$ is representation-infinite. Therefore, $B$ is isomorphic to either $\mathbf{A}_2^{(+)}$ or $N(n)$ with $n\ge 3$ if $\mathbf{A}_3^{\varepsilon}\otimes B$ is representation-finite.

By Corollary \ref{cor::A2-An}, $\mathbf{A}_3^\varepsilon\otimes \mathbf{A}_2^{(+)}$ is representation-finite for any choice of $\varepsilon$.
We then show that $\mathbf{A}_3^{\varepsilon}\otimes N(n)$ is also representation-finite for each $\varepsilon$. Up to opposite algebras, it suffices to consider two cases: $(++)$ and $(+-)$. In fact, it follows from \cite[Theorem 6.1]{LS-tame-tensor-product} that $\mathbf{A}_3^{(++)}\otimes N(n)$ is representation-finite. The tensor product $\mathbf{A}_3^{(+-)}\otimes N(n)$ is presented by
$$
\vcenter{\xymatrix@C=1.5cm@R=0.7cm{
\circ \ar[r]^a\ar[d]^b\ar@{.}[dr]&\circ\ar[r]^a\ar[d]^b\ar@{.}[dr]&\circ\ar[r]^a\ar[d]^b\ar@{.}[dr]&\circ\ar[r]^a\ar[d]^b\ar@{.}[dr]&\cdots\ar[r]^a\ar@{.}[dr]&\circ\ar[r]^a\ar[d]^b\ar@{.}[dr]&\circ\ar[d]^b\\
\circ \ar[r]^a  &\circ\ar[r]^a  &\circ\ar[r]^a&\circ\ar[r]^a & \cdots\ar[r]^a &\circ\ar[r]^a &\circ  \\
\circ\ar[r]^a\ar[u]_b\ar@{.}[ur]&\circ  \ar[r]^a\ar[u]_b\ar@{.}[ur]&\circ\ar[r]^a \ar[u]_b\ar@{.}[ur]&\circ\ar[r]^a\ar[u]_b\ar@{.}[ur]&\cdots\ar[r]^a&\circ \ar[u]_b\ar@{.}[ur]\ar[r]^a&\circ\ar[u]_b
}}.
$$
with $a^2=0$ and $ab=ba$. We observe that $\mathbf{A}_3^{(+-)}\otimes N(n)$ does not contain any tame concealed algebra as a quotient. In fact, if $\mathbf{A}_3^{(+-)}\otimes N(n)$ has a tame concealed algebra (say, $C$) of type $\widetilde{\mathbb{D}}_n$ or $\widetilde{\mathbb{E}}_n$ as a quotient algebra, then $C$ can be embedded in a truncation of $\mathbf{A}_3^{(+-)}\otimes N(n)$ consisting of 3 rows and 8 columns in its quiver, see Remark \ref{rem::shape-4*8}. (Here, we assume that $n$ is large enough to take truncation. If not so, then the case is easy to check.) Such a $3\times 8$ truncation of $\mathbf{A}_3^{(+-)}\otimes N(n)$ is always isomorphic to $\mathbf{A}_3^{(+-)}\otimes N(8)$, and it is not difficult to check that $\mathbf{A}_3^{(+-)}\otimes N(8)$ contains no tame concealed algebras. Thus, $\mathbf{A}_3^{(+-)}\otimes N(n)$ is representation-finite following Proposition \ref{prop::Tits-form}.
\end{proof}

\begin{remark}\label{rem::GAP-code}
One may use the software GAP 4.10.2\footnote{Some bugs have appeared in the latest version of GAP, i.e., GAP 4.12.2, so that the command $\mathsf{IsWeaklyPositiveUnitForm(-)}$ can only be used in GAP 4.10.2 or earlier versions.} to find that the Tits form of $\mathbf{A}_3^{(+-)}\otimes N(8)$ is weakly positive. The code we use is as follows.
\begin{quote}
\begin{lstlisting}[language=Python]
LoadPackage("qpa");;
# $A \simeq \mathbf{A}_3^{(+-)}$
Q_A:=Quiver(["v1","v2","v3"], [["v1","v2","x"],["v3","v2","y"]]);;
KQ_A:=PathAlgebra(Rationals,Q_A);;
AssignGeneratorVariables(KQ_A);;
A:=KQ_A;;
# $B \simeq N(8)$
Q_B:=Quiver(["v1","v2","v3","v4","v5","v6","v7","v8"], [["v1","v2","a"],["v2","v3","b"],["v3","v4","c"],["v4","v5","d"],["v5","v6","e"],["v6","v7","f"],["v7","v8","g"]]);;
KQ_B:=PathAlgebra(Rationals,Q_B);;
AssignGeneratorVariables(KQ_B);;
rel_B:=[a*b,b*c,c*d,d*e,e*f,f*g];;
B:=KQ_B/rel_B;;
# $T=\mathbf{A}_3^{(+-)} \otimes N(8)$
T:=TensorProductOfAlgebras(A,B);;
q_T:=TitsUnitFormOfAlgebra(T);;
IsWeaklyPositiveUnitForm(q_T);
# Output
true
\end{lstlisting}
\end{quote}
\end{remark}

\begin{lemma}\label{lem::A4}
Let $B$ be a non-local representation-finite simply connected algebra. Then, $\mathbf{A}_4^{\varepsilon}\otimes B$ is representation-finite if and only if $B$ is isomorphic to $\mathbf{A}_2^{(+)}$.
\end{lemma}
\begin{proof}
If $B\simeq \mathbf{A}_2^{(+)}$, then $\mathbf{A}_4^{\varepsilon}\otimes B$ is representation-finite by Corollary \ref{cor::A2-An}. We then show that $\mathbf{A}_4^{\varepsilon}\otimes N(n)$ with $n\ge 3$ is representation-infinite, and the statement follows from Lemma \ref{lem::A3}. It only needs to consider 4 cases: $(+++)$, $(+-+)$, $(-++)$ and $(++-)$.
\begin{itemize}
\item Set $\varepsilon=(+++)$. It follows from \cite[Theorem 6.1]{LS-tame-tensor-product} that $\mathbf{A}_4^{(+++)}\otimes N(3)$ is representation-infinite.

\item Set $\varepsilon=(+-+)$. Then, $\mathbf{A}_4^{(+-+)}\otimes N(3)$ is presented by the quiver
$$
\vcenter{\xymatrix@C=1.5cm@R=0.7cm{
\circ \ar[r]^a\ar[d]^b\ar@{.}[dr]&\bullet\ar[d]^b&\circ\ar[l]_a\ar[r]^a\ar[d]^b\ar@{.}[dr]\ar@{.}[dl]&\circ\ar[d]^b\\
\bullet \ar[r]^a\ar[d]^b\ar@{.}[dr]&\bullet\ar[d]^b &\bullet\ar[l]_a\ar[r]^a\ar[d]^b\ar@{.}[dl]\ar@{.}[dr]&\bullet\ar[d]^b\\
\circ\ar[r]^a&\circ&\bullet\ar[r]^a\ar[l]_a&\circ
}}
$$
with $b^2=0$ and $ab=ba$. Here, the black points indicate a tame concealed algebra of type $\widetilde{\mathbb{D}}_5$ such that $\mathbf{A}_4^{(+-+)}\otimes N(3)$ is representation-infinite.

\item Set $\varepsilon=(-++)$. Then, $\mathbf{A}_4^{(-++)}\otimes N(3)$ is presented by the quiver
$$
\vcenter{\xymatrix@C=1.5cm@R=0.7cm{
\bullet\ar[d]^b&\circ\ar[l]_a\ar[r]^a\ar[d]^b\ar@{.}[dl]\ar@{.}[dr]&\circ\ar[r]^a\ar[d]^b\ar@{.}[dr]&\bullet\ar[d]^b\\
\bullet \ar[d]^b &\bullet\ar[l]_a\ar[r]^a\ar[d]^b\ar@{.}[dr]\ar@{.}[dl]&\bullet\ar[r]^a\ar[d]^b\ar@{.}[dr]&\bullet\ar[d]^b\\
\circ &\bullet\ar[l]_a\ar[r]^a&\bullet\ar[r]^a&\circ
}}
$$
with $b^2=0$ and $ab=ba$. The black points indicate a tame concealed algebra of type $\widetilde{\mathbb{E}}_7$ (see No.13 in \cite{B-critical}) such that $\mathbf{A}_4^{(-++)}\otimes N(3)$ is representation-infinite.

\item Set $\varepsilon=(++-)$. Then, $\mathbf{A}_4^{(++-)}\otimes N(3)$ is presented by the quiver
$$
\vcenter{\xymatrix@C=1.5cm@R=0.7cm{
\circ \ar[r]^a\ar[d]^b\ar@{.}[dr]&\circ\ar[r]^a\ar[d]^b\ar@{.}[dr]&\bullet \ar[d]^b&\circ\ar[d]^b\ar[l]_a\ar@{.}[dl]\\
\bullet \ar[r]^a\ar[d]^b\ar@{.}[dr]&\bullet\ar[r]^a\ar[d]^b\ar@{.}[dr]&\bullet \ar[d]^b&\bullet\ar[d]^b\ar[l]_a\ar@{.}[dl]\\
\bullet\ar[r]^a&\bullet\ar[r]^a&\circ&\circ\ar[l]_a
}}
$$
with $b^2=0$ and $ab=ba$. The black points in the quiver indicate a tame concealed algebra of type $\widetilde{\mathbb{D}}_6$ such that $\mathbf{A}_4^{(++-)}\otimes N(3)$ is representation-infinite.
\end{itemize}
By taking the surjection $N(n)\rightarrow N(3)$, we find that $\mathbf{A}_4^{\varepsilon}\otimes N(n)$ with $n\ge 3$ is representation-infinite in all the 4 cases.
\end{proof}

\begin{proposition}
Let $B$ be a non-local representation-finite simply connected algebra.  Then, $\mathbf{A}_n^{\varepsilon}\otimes B$ with $n\ge 5$ is representation-infinite for any orientation $\varepsilon$.
\end{proposition}
\begin{proof}
There exists a surjection $\mathbf{A}_n^{\varepsilon}\otimes B\to \mathbf{A}_5^{\varepsilon}\otimes \mathbf{A}_2^{(+)}$, and the latter one is representation-infinite by Corollary \ref{cor::A2-An}.
\end{proof}

\subsection{Path algebra of type $D$}
Let $\mathbb{D}_n$ ($n\geq 4$) be the Dynkin diagram of type $D$:
$$
\vcenter{\xymatrix@C=1cm@R=0.1cm{
1\ar@{-}[dr]&&&&&\\
&3\ar@{-}[r]&4\ar@{-}[r]&\cdots\ar@{-}[r]& n-1\ar@{-}[r]& n\\
2\ar@{-}[ur]&&&&&}}.
$$
The orientation of $\mathbb{D}_n$ is defined as $\varepsilon:=(\varepsilon_1,\varepsilon_2,\ldots, \varepsilon_{n-1})$ with
$$
\left\{\begin{array}{ll}
\varepsilon _1=+ & \text{if } 1 \longrightarrow 3 , \\
\varepsilon _1=- & \text{if } 3 \longrightarrow 1 .
\end{array}\right.
\quad
\left\{\begin{array}{ll}
\varepsilon _2=+ & \text{if } 2 \longrightarrow 3 , \\
\varepsilon _2=- & \text{if } 3 \longrightarrow 2 .
\end{array}\right.
\quad
\left\{\begin{array}{ll}
\varepsilon _i=+ & \text{if } i \longrightarrow i+1 , \\
\varepsilon _i=- & \text{if } i+1 \longrightarrow i .
\end{array}\right. (i\ge 3)
$$
We denote by $\mathbf{D}_n^\varepsilon$ the path algebra of type $D$ associated with the orientation $\varepsilon$.

\begin{proposition}\label{D-case}
Let $B$ be a non-local representation-finite simply connected algebra. Then, the tensor product $\mathbf{D}_n^{\varepsilon}\otimes B$ is representation-infinite for any $n\geq 4$.
\end{proposition}
\begin{proof}
It is enough to show that $\mathbf{D}_4^{\varepsilon}\otimes\mathbf{A}_2^{(+)}$ is representation-infinite. In fact, this follows from Proposition \ref{prop::A2} since $\mathbb{D}_4$ is included in the family (IT) in \cite{LS-tame-triangular-matrix}.
\end{proof}

Let $\mathbb{E}_6$, $\mathbb{E}_7$ and $\mathbb{E}_8$ be the Dynkin diagrams of type $E$. These contain $\mathbb{D}_4$ as a subgraph, we have the following corollary from Proposition \ref{D-case}.

\begin{corollary}
Let $B$ be a non-local representation-finite simply connected algebra and $\mathbf{E}$ a path algebra of type $E$. Then, the tensor product $\mathbf{E}\otimes B$ is representation-infinite.
\end{corollary}

We have also obtained a slight generalization of \cite[Theorem 4]{Aihara-Honma}, in which the authors deal with the representation-finiteness of $T_n(B)$.

\begin{corollary}\label{cor::path-alg-case}
Suppose $A$ is a representation-finite hereditary algebra with $n\geq 2$ simple modules. Then, the following statements hold.
\begin{enumerate}
\item[\rm{(1)}] Let $B$ be a non-local hereditary algebra. Then $A\otimes B$ is representation-finite if and only if $A\simeq \mathbf{A}_2^{(+)}$ and $B$ is isomorphic to one of path algebras of $\mathbb{A}_2$, $\mathbb{A}_3$ and $\mathbb{A}_4$.

\item[\rm{(2)}] Let $n\geq 3$ and $B$ be non-local simply connected but not hereditary algebra. Then, $A\otimes B$ is representation-finite if and only if $A$ is isomorphic to a path algebra of $\mathbb{A}_3$ and $B\simeq N(m)$ for some $m\ge 3$.
\end{enumerate}
\end{corollary}

\subsection{Result}
We have the following result related to our classification.
\begin{theorem}\label{theo::not-nakayama}
Let $A$ and $B$ be two non-local simply connected algebras. Then, the following statements hold.
\begin{enumerate}
\item[\rm{(1)}] If both $A$ and $B$ are not Nakayama algebras, then $A\otimes B$ is representation-infinite.

\item[\rm{(2)}] If $A$ is a Nakayama algebra with $\mathsf{rad}^2\ A\neq 0$ and $B$ is not a Nakayama algebra, then $A\otimes B$ is representation-infinite.

\item[\rm{(3)}] If both $A$ and $B$ are Nakayama algebras with radical square not zero, then $A\otimes B$ is representation-infinite.
\end{enumerate}
\end{theorem}
\begin{proof}
In all three cases, we notice that there is a surjection $A\otimes B\to \mathbf{A}_3^\varepsilon\otimes \mathbf{A}_3^\omega$ for some choices of $\varepsilon$ and $\omega$. Then, the statements follow from Lemma \ref{lem::A3-A3}.
\end{proof}

\section{Nakayama algebras}
We have reduced the general problem on $A\otimes B$ to the cases that at least one of $A$ and $B$ is a Nakayama algebra with radical square zero, see Theorem \ref{theo::not-nakayama}. Up to isomorphism, we assume in this section that $A\simeq N(n)$ with $n\ge 3$. Moreover, we assume $|B|\ge 3$ since $B\simeq \mathbf{A}_2^{(+)}$ if $|B|=2$.

Let us start with the following easy observation.
\begin{proposition}
The tensor product $N(n)\otimes N(m)$ with $n, m\ge 3$ is representation-finite.
\end{proposition}
\begin{proof}
By constructing the quiver presentation of $N(n)\otimes N(m)$, we find that it is a special biserial algebra. It is then shown in \cite{ABM} that $N(n)\otimes N(m)$ is representation-finite.
\end{proof}

In the following, we divide the problem on $N(n)\otimes B$ into two disjoint cases: $B$ is a Nakayama algebra or $B$ is not a Nakayama algebra.

\subsection{$B$ is Nakayama}
We assume that $B$ is a simply connected Nakayama algebra and is not radical square zero. If $|B|=3$, then $B\simeq \mathbf{A}_3^{(++)}$ and this case has been determined in Corollary \ref{cor::path-alg-case}, i.e., $N(n)\otimes B$ is representation-finite for any $n\ge 3$. Suppose $|B|=m\ge 4$ in this subsection. We mention that the quiver of $N(n)\otimes B$ is displayed as
$$
\Delta_{n,m}:\vcenter{\xymatrix@C=1.5cm@R=0.85cm{
\circ \ar[r]^a\ar[d]^b\ar@{.}[dr]&\circ\ar[r]^a\ar[d]^b\ar@{.}[dr]&\circ\ar[r]^a\ar[d]^b\ar@{.}[dr]&\circ\ar[r]^a\ar[d]^b\ar@{.}[dr]&\cdots\ar[r]^a\ar@{.}[dr]&\circ\ar[r]^a\ar[d]^b\ar@{.}[dr]&\circ\ar[d]^b\\
\circ \ar[r]^a\ar[d]^b\ar@{.}[dr]&\circ\ar[r]^a\ar[d]^b\ar@{.}[dr]&\circ\ar[r]^a\ar[d]^b\ar@{.}[dr]&\circ\ar[r]^a\ar[d]^b\ar@{.}[dr]& \cdots\ar[r]^a\ar@{.}[dr]&\circ\ar[r]^a\ar[d]^b\ar@{.}[dr]&\circ \ar[d]^b\\
\circ \ar[r]^a\ar[d]^b\ar@{.}[dr]&\circ\ar[r]^a\ar[d]^b\ar@{.}[dr]&\circ\ar[r]^a\ar[d]^b\ar@{.}[dr]&\circ\ar[r]^a\ar[d]^b\ar@{.}[dr]& \cdots\ar[r]^a\ar@{.}[dr]&\circ\ar[r]^a\ar[d]^b\ar@{.}[dr]&\circ \ar[d]^b\\
\vdots \ar[d]^b\ar@{.}[dr]&\vdots \ar[d]^b\ar@{.}[dr]&\vdots \ar[d]^b\ar@{.}[dr]&\vdots \ar[d]^b\ar@{.}[dr]& \cdots\ar@{.}[dr]& \vdots \ar[d]^b\ar@{.}[dr]&\vdots \ar[d]^b\\
\circ\ar[r]^a&\circ  \ar[r]^a &\circ\ar[r]^a &\circ\ar[r]^a&\cdots\ar[r]^a&\circ \ar[r]^a&\circ
}}.
$$

\begin{lemma}\label{lem::N3-A4}
If $\mathsf{rad}^3\ B\neq 0$, then $N(n)\otimes B$ is representation-infinite for any $n\ge 3$.
\end{lemma}
\begin{proof}
We show that $N(3)\otimes B$ is representation-infinite. Since $\mathsf{rad}^3\ B\neq 0$, $N(3)\otimes B$ must admit a quotient algebra with the following quiver
$$
\vcenter{\xymatrix@C=1.5cm@R=0.7cm{
\circ \ar[r]^a\ar[d]^b\ar@{.}[dr]&\circ\ar[r]^a\ar[d]^b\ar@{.}[dr]&\bullet\ar[r]^a\ar[d]^b\ar@{.}[dr]&\bullet\ar[d]^b\\
\bullet \ar[r]^a\ar[d]^b\ar@{.}[dr]&\bullet\ar[r]^a\ar[d]^b\ar@{.}[dr]&\bullet\ar[r]^a\ar[d]^b\ar@{.}[dr]&\bullet\ar[d]^b\\
\bullet\ar[r]^a&\bullet\ar[r]^a&\circ\ar[r]^a&\circ
}}
$$
and $a^3\neq0, b^2=0$, $ab=ba$. The black points in the quiver indicate a tame concealed algebra of type $\widetilde{\mathbb{D}}_7$ such that the quotient of $N(3)\otimes B$ is representation-infinite.
\end{proof}

If $\mathsf{rad}^3\ B=0$, then $N(n)\otimes B$ is a quotient algebra of $\Lambda_{n,m}:=\mathsf{k}\Delta_{n,m}/\langle a^3,b^2,ab-ba\rangle$, for some $n\ge 3, m\ge 4$.
Since $a^3=0$ and $b^2=0$, it is easy to find that $\Lambda_{n,m}$ cannot have a tame concealed algebra of type $\widetilde{\mathbb{D}}_n$ as a quotient algebra. We have the following results.

\begin{lemma}\label{lem::N3-B1}
If $B$ has the bound quiver algebra
\begin{equation}
\mathbf{B}_1:=\mathsf{k}\left(
\xymatrix@C=1cm@R=0.5cm{
1\ar[r]^\alpha&2\ar[r]^\beta&3\ar[r]^\gamma&4\ar[r]^\delta&5}
\right)
\bigg/
\langle \alpha\beta\gamma, \beta\gamma\delta \rangle
\end{equation}
as a quotient, then $N(3)\otimes B$ is representation-infinite.
\end{lemma}
\begin{proof}
The tensor product $N(3)\otimes \mathbf{B}_1$ is given by the quiver
$$
\vcenter{\xymatrix@C=1.5cm@R=0.7cm{
\circ \ar[r]^a\ar[d]^b\ar@{.}[dr]&\circ\ar@{.}[dr]\ar[r]^a\ar[d]^b &\bullet\ar[d]^b\ar[r]^a\ar@{.}[dr]&\bullet \ar[r]^a\ar@{.}[dr]\ar[d]^b &\bullet \ar[d]^b \\
\circ \ar[r]^a\ar[d]^b\ar@{.}[dr]&\bullet \ar@{.}[dr]\ar[r]^a\ar[d]^b &\bullet\ar[d]^b\ar[r]^a\ar@{.}[dr]&\bullet\ar[r]^a\ar[d]^b\ar@{.}[dr] &\circ\ar[d]^b \\
\bullet\ar[r]^a&\bullet \ar[r]^a&\bullet \ar[r]^a&\circ\ar[r]^a&\circ
}}
$$
with $a^3=b^2=0$ and $ab=ba$. The black points give a tame concealed algebra of type $\widetilde{\mathbb{E}}_8$ (see No.99 in \cite{B-critical}). This implies that $N(3)\otimes B$ is representation-infinite.
\end{proof}

\begin{lemma}\label{lem::N3-B2}
If $B$ has the bound quiver algebra
\begin{equation}
\mathbf{B}_2:=\mathsf{k}\left(
\xymatrix@C=0.7cm@R=0.5cm{
1\ar[r]^\alpha&2\ar[r]^\beta&3\ar[r]^\gamma&4\ar[r]^\delta&5\ar[r]^\sigma &6}
\right)
\bigg/
\langle \alpha\beta\gamma, \gamma\delta \rangle
\end{equation}
as a quotient, then $N(3)\otimes B$ is representation-infinite.
\end{lemma}
\begin{proof}
The tensor product $N(3)\otimes \mathbf{B}_2$ is given by the quiver
$$
\vcenter{\xymatrix@C=1.5cm@R=0.7cm{
\circ \ar[r]^a\ar[d]^b\ar@{.}[dr]&\circ\ar@{.}[dr]\ar[r]^a\ar[d]^b &\circ\ar[d]^b\ar[r]^c\ar@{.}[dr]&\bullet \ar[r]^c\ar@{.}[dr]\ar[d]^b &\bullet \ar[r]^d\ar@{.}[dr]\ar[d]^b &\bullet \ar[d]^b \\
\circ \ar[r]^a\ar[d]^b\ar@{.}[dr]&\bullet \ar@{.}[dr]\ar[r]^a\ar[d]^b &\bullet\ar[d]^b\ar[r]^c\ar@{.}[dr]&\bullet\ar[r]^c\ar[d]^b\ar@{.}[dr]&\circ \ar[r]^d\ar@{.}[dr]\ar[d]^b &\circ\ar[d]^b \\
\bullet\ar[r]^a&\bullet \ar[r]^a&\bullet \ar[r]^c&\circ\ar[r]^c&\circ\ar[r]^d&\circ
}}
$$
with $a^2c=b^2=c^2=0$, $ab=ba$, $bc=cb$ and $bd=db$. The black points give a tame concealed algebra of type $\widetilde{\mathbb{E}}_8$ (see No.17 in \cite{B-critical}). This implies that $N(3)\otimes B$ is representation-infinite.
\end{proof}

\begin{proposition}\label{prop::N3-Nakayama}
Suppose $\mathsf{rad}^3\ B=0$. Then, $N(3)\otimes B$ is representation-finite if and only if $B$ does not contain one of $\mathbf{B}_1$, $\mathbf{B}_2$ and their opposite algebras as a quotient algebra.
\end{proposition}
\begin{proof}
Suppose $B$ does not contain one of $\mathbf{B}_1$, $\mathbf{B}_2$, $\mathbf{B}_1^{\mathsf{op}}$, $\mathbf{B}_2^{\mathsf{op}}$ as a quotient algebra. If $N(3)\otimes B$ has a tame concealed algebra (say, $C$) of type $\widetilde{\mathbb{E}}_n$ as a quotient algebra, then $C$ can be embedded in a truncation of $N(3)\otimes B$ consisting of 3 rows and 8 columns in $\Delta_{3,m}$, see Remark \ref{rem::shape-4*8}. We assume that $m$ is large enough to take truncation. If not so, then the case is some quotient of the cases with larger $m$.

Set $|B|=8$. All possible choices of $B$ are listed as follows.
\begin{enumerate}
\item $\xymatrix@C=0.5cm{1\ar[r]^{a_1}&2\ar[r]^{a_2}&3\ar[r]^{a_3}&4\ar[r]^{a_4}&5\ar[r]^{a_5}&6\ar[r]^{a_6}&7\ar[r]^{a_7}&8}$ with $a_1a_2a_3=a_3a_4=a_4a_5=a_5a_6a_7=0$.

\item [$(1^1)$]
$\xymatrix@C=0.5cm{1\ar[r]^{a_1}&2\ar[r]^{a_2}&3\ar[r]^{a_3}&4\ar[r]^{a_4}&5\ar[r]^{a_5}&6\ar[r]^{a_6}&7\ar[r]^{a_7}&8}$ with $a_1a_2a_3=a_3a_4=a_4a_5=a_6a_7=0$.

\item [$(1^2)$]
$\xymatrix@C=0.5cm{1\ar[r]^{a_1}&2\ar[r]^{a_2}&3\ar[r]^{a_3}&4\ar[r]^{a_4}&5\ar[r]^{a_5}&6\ar[r]^{a_6}&7\ar[r]^{a_7}&8}$ with $a_1a_2a_3=a_3a_4=a_4a_5=a_5a_6=0$.

\item [$(1^3)$]
$\xymatrix@C=0.5cm{1\ar[r]^{a_1}&2\ar[r]^{a_2}&3\ar[r]^{a_3}&4\ar[r]^{a_4}&5\ar[r]^{a_5}&6\ar[r]^{a_6}&7\ar[r]^{a_7}&8}$ with $a_1a_2=a_3a_4=a_4a_5=a_5a_6a_7=0$.

\item [$(1^4)$]
$\xymatrix@C=0.5cm{1\ar[r]^{a_1}&2\ar[r]^{a_2}&3\ar[r]^{a_3}&4\ar[r]^{a_4}&5\ar[r]^{a_5}&6\ar[r]^{a_6}&7\ar[r]^{a_7}&8}$ with $a_2a_3=a_3a_4=a_4a_5=a_5a_6a_7=0$.

\item 
$\xymatrix@C=0.5cm{1\ar[r]^{a_1}&2\ar[r]^{a_2}&3\ar[r]^{a_3}&4\ar[r]^{a_4}&5\ar[r]^{a_5}&6\ar[r]^{a_6}&7\ar[r]^{a_7}&8}$ with $a_1a_2=a_2a_3a_4=a_4a_5=a_5a_6=0$.

\item [$(2^{\mathsf{op}})$]
$\xymatrix@C=0.5cm{1\ar[r]^{a_1}&2\ar[r]^{a_2}&3\ar[r]^{a_3}&4\ar[r]^{a_4}&5\ar[r]^{a_5}&6\ar[r]^{a_6}&7\ar[r]^{a_7}&8}$ with $a_2a_3=a_3a_4=a_4a_5a_6=a_6a_7=0$.

\item $\xymatrix@C=0.5cm{1\ar[r]^{a_1}&2\ar[r]^{a_2}&3\ar[r]^{a_3}&4\ar[r]^{a_4}&5\ar[r]^{a_5}&6\ar[r]^{a_6}&7\ar[r]^{a_7}&8}$ with $a_1a_2=a_3a_4=a_5a_6=0$.

\item [$(3^{\mathsf{op}})$]
$\xymatrix@C=0.5cm{1\ar[r]^{a_1}&2\ar[r]^{a_2}&3\ar[r]^{a_3}&4\ar[r]^{a_4}&5\ar[r]^{a_5}&6\ar[r]^{a_6}&7\ar[r]^{a_7}&8}$ with $a_2a_3=a_4a_5=a_6a_7=0$.

\item $\xymatrix@C=0.5cm{1\ar[r]^{a_1}&2\ar[r]^{a_2}&3\ar[r]^{a_3}&4\ar[r]^{a_4}&5\ar[r]^{a_5}&6\ar[r]^{a_6}&7\ar[r]^{a_7}&8}$ with $a_2a_3=a_4a_5=a_5a_6=0$.

\item [$(4^{\mathsf{op}})$]
$\xymatrix@C=0.5cm{1\ar[r]^{a_1}&2\ar[r]^{a_2}&3\ar[r]^{a_3}&4\ar[r]^{a_4}&5\ar[r]^{a_5}&6\ar[r]^{a_6}&7\ar[r]^{a_7}&8}$ with $a_2a_3=a_3a_4=a_5a_6=0$.

\item $\xymatrix@C=0.5cm{1\ar[r]^{a_1}&2\ar[r]^{a_2}&3\ar[r]^{a_3}&4\ar[r]^{a_4}&5\ar[r]^{a_5}&6\ar[r]^{a_6}&7\ar[r]^{a_7}&8}$ \\ 
with $a_ia_{i+1}=0$ for all but one $i\in \{1,2,\cdots, 6\}$.

\item $\xymatrix@C=0.5cm{1\ar[r]^{a_1}&2\ar[r]^{a_2}&3\ar[r]^{a_3}&4\ar[r]^{a_4}&5\ar[r]^{a_5}&6\ar[r]^{a_6}&7\ar[r]^{a_7}&8}$ \\ 
with $a_ia_{i+1}=0$ for all but two $i\neq j\in \{1,2,\cdots, 6\}$.
\end{enumerate}
Here, the cases $(1^i)$'s are quotients of $(1)$ and the case $(i^{\mathsf{op}})$ is the opposite algebra of $(i)$. Moreover, the cases $(5)$ and $(6)$ are quotients of some cases in $(1)-(4)$. Similar to Remark \ref{rem::GAP-code}, one may use GAP 4.10.2
to check that the Tits form of $N(3)\otimes B$ is weakly positive in each case of $(1)-(4)$. 

We observe that each $3\times 8$ truncation of $N(3)\otimes B$ is isomorphic to one of the algebras listed above. It turns out that $N(3)\otimes B$ is representation-finite by Proposition \ref{prop::Tits-form}.
\end{proof}

\begin{lemma}\label{lem::N4-B3}
If $B$ has the bound quiver algebra
\begin{equation}\label{B_0}
\mathbf{B}_3:=\mathsf{k}\left(
\xymatrix@C=1cm@R=0.5cm{
1\ar[r]^\alpha&2\ar[r]^\beta&3\ar[r]^\gamma&4}
\right)
\bigg/
\langle\alpha\beta\gamma\rangle
\end{equation}
as a quotient, then $N(n)\otimes B$ is representation-infinite for any $n\geq 4$.
\end{lemma}
\begin{proof}
We show that $N(4)\otimes \mathbf{B}_3$ is representation-infinite. In fact, $N(4)\otimes \mathbf{B}_3$ is isomorphic to $\Lambda_{4,4}$. The black points in $\Delta_{4,4}$, i.e.,
$$
\vcenter{\xymatrix@C=1.5cm@R=0.7cm{
\circ \ar[r]^a\ar[d]^b\ar@{.}[dr]&\circ\ar[r]^a\ar[d]^b\ar@{.}[dr]&\circ\ar[r]^a\ar[d]^b\ar@{.}[dr]&\bullet\ar[d]^b\\
\circ \ar[r]^a\ar[d]^b\ar@{.}[dr]&\bullet\ar[r]^a\ar[d]^b\ar@{.}[dr]&\bullet\ar[r]^a\ar[d]^b\ar@{.}[dr]&\bullet\ar[d]^b\\
\bullet\ar[r]^a\ar[d]^b\ar@{.}[dr]&\bullet\ar[r]^a\ar[d]^b\ar@{.}[dr]&\bullet\ar[r]^a\ar[d]^b\ar@{.}[dr]&\circ\ar[d]^b\\
\bullet\ar[r]^a&\circ\ar[r]^a&\circ\ar[r]^a&\circ
}}
$$
give a tame concealed algebra of type $\widetilde{\mathbb{E}}_7$ (see No.12 in \cite{B-critical}). Hence,  $N(4)\otimes \mathbf{B}_3$ is representation-infinite.
\end{proof}

\begin{proposition}
Suppose $\mathsf{rad}^3\ B=0$. Then, $N(n)\otimes B$ with $n\ge 4$ is representation-finite if and only if $B$ does not contain $\mathbf{B}_3$ as a quotient algebra.
\end{proposition}
\begin{proof}
We only show the sufficiency. If $B$ does not contain $\mathbf{B}_3$ as a quotient algebra, then we have $\alpha\beta=0$ or $\beta\gamma=0$ for each path
$\xymatrix@C=0.6cm@R=0.5cm{\circ\ar[r]^\alpha&\circ\ar[r]^\beta&\circ\ar[r]^\gamma&\circ}$ in $B$. It suffices to consider $B$ as a quotient of $\Gamma_m$, where $\Gamma_m$ ($m\ge 4$) is given by
$$
\xymatrix@C=1cm@R=0.5cm{1\ar[r]^{a_1}&2\ar[r]^{a_2}&3\ar[r]^{a_1}&4\ar[r]^{a_2}&5\ar[r]^{a_1}&6\ar[r]^{a_2}&\cdots \ar[r] &m-1\ar[r] &m}
$$
with $a_1a_2=0$ and $a_2a_1\neq 0$. We show that $N(n)\otimes \Gamma_m$ is representation-finite.

Since $N(n)\otimes \Gamma_m$ is a quotient of $\Lambda_{n,m}$, it also cannot have a tame concealed algebra of type $\widetilde{\mathbb{D}}_n$ as a quotient algebra. Similar to the situation of Lemma \ref{lem::A3}, if $N(n)\otimes \Gamma_m$ admits a tame concealed algebra of type $\widetilde{\mathbb{E}}_n$ as a quotient, then the quotient must be embedded in a $4\times 8$ or $8\times 4$ truncation of $N(n)\otimes \Gamma_m$, see Remark \ref{rem::shape-4*8}. (Here, we assume that $n$ and $m$ are large enough to take truncation.) Up to isomorphism, we only need to consider $N(4)\otimes \Gamma_8$ and $N(8)\otimes \Gamma_4$, and both cases do not contain a tame concealed algebra of type $\widetilde{\mathbb{E}}_n$ as a quotient. (One may also use the software GAP 4.10.2 to find that the Tits forms of $N(4)\otimes \Gamma_8$ and $N(8)\otimes \Gamma_4$ are weakly positive, see Remark \ref{rem::GAP-code}.) Thus, $N(n)\otimes \Gamma_m$ is representation-finite by Proposition \ref{prop::Tits-form}.
\end{proof}

In summary, we have the following result related to our classification.
\begin{theorem}\label{Theo::Nn-naka-rad-not-zero}
Let $A\simeq N(n)$ with $n\ge 3$ and $B$ a Nakayama algebra with $|B|\ge 3$ and $\mathsf{rad}^2\ B\neq 0$. Then, $A\otimes B$ is representation-finite if and only if one of the following holds:
\begin{enumerate}
\item[\rm{(1)}] $n=3$, $\mathsf{rad}^3\ B=0$ and $B$ does not contain one of $\mathbf{B}_1$, $\mathbf{B}_2$ and their opposite algebras as a quotient algebra.
\item[\rm{(2)}] $n\ge 4$ and $B$ does not contain $\mathbf{B}_3$ as a quotient algebra.
\end{enumerate}
\end{theorem}

\subsection{$B$ is not Nakayama}
In this subsection, we consider the case that $B$ is a simply connected algebra, but is not a Nakayama algebra. Then, the quiver $Q_B$ of $B$ contains $\xymatrix@C=0.6cm{\circ&\circ\ar[l]\ar[r]&\circ}$ or $\xymatrix@C=0.6cm{\circ\ar[r]&\circ &\circ\ar[l]}$ as a subquiver.

\begin{proposition}
If $|B|=3$, then $N(n)\otimes B$ is representation-finite for any $n\geq 3$.
\end{proposition}
\begin{proof}
Let $C=\mathsf{k}Q_C/\mathcal{I}$ be a bound quiver algebra.
If $Q_C$ is an acyclic quiver whose underlying graph is a triangle, that is,
$$
\vcenter{\xymatrix@C=0.3cm@R=0.4cm{&\circ\ar@{-}[ld]\ar@{-}[rd]&\\\circ\ar@{-}[rr]&&\circ}},
$$
then the fundamental group of $C$ is isomorphic to $\mathbb{Z}$ such that $C$ is not simply connected.
Since $B$ is simply connected but not Nakayama and $|B|=3$, we have $B\simeq \mathbf{A}_3^{(+-)}$ or $B\simeq \mathbf{A}_3^{(-+)}$. Then, we obtain the statement from Corollary \ref{cor::path-alg-case}.
\end{proof}

\begin{proposition}\label{prop::N3-B4}
If the underlying graph of the quiver $Q_B$ of $B$ has
$$
\vcenter{\xymatrix@C=0.5cm@R=0.01cm{
\circ\ar@{-}[dr]&&\\
&\circ\ar@{-}[r]&\circ\\
\circ\ar@{-}[ur]&&}}
$$
as a subgraph, then $N(n)\otimes B$ is representation-infinite for any $n\geq 3$.
\end{proposition}
\begin{proof}
It is enough to check the case $|B|=4$. If $B\simeq \mathbf{D}_4^{\varepsilon}$ for some $\varepsilon$, then $N(n)\otimes B$ is representation-infinite by Corollary \ref{cor::path-alg-case}. Otherwise, $B$ has either $\mathbf{B}_4$ or $\mathbf{B}_4^{\mathsf{op}}$ as a quotient algebra, where $\mathbf{B}_4$ is presented as
$$
\mathbf{B}_4:=\mathsf{k}\left(
\vcenter{\xymatrix@C=0.5cm@R=0.1cm{
1\ar[dr]^\alpha&&\\
&3\ar[dl]^\beta&4\ar[l]_\gamma\\
2&&}}
\right) \bigg/ \langle \alpha\beta, \gamma\beta\rangle.
$$
It suffices to show that $N(3)\otimes \mathbf{B}_4$ is representation-infinite. In fact, $N(3)\otimes \mathbf{B}_4$ is given by
$$
\vcenter{\xymatrix@C=0.7cm@R=0.3cm{
\circ \ar[rr]^a\ar[dd]^b&&\bullet\ar[dl]_c\ar[dd]^b&\circ\ar[l]_a\ar[dd]^b\\
&\bullet\ar@/_0.3cm/[dd]_b&&\\
\bullet \ar@{..>}[rr]^a\ar[dd]^b&&\bullet\ar[dl]_c\ar[dd]^b&\bullet\ar[l]_a\ar[dd]^b\\
&\circ\ar@/_0.3cm/[dd]_b&&\\
\bullet \ar@{..>}[rr]^a&&\circ\ar[dl]_c&\bullet\ar[l]_a\\
&\circ&&
}}
$$
with $ac=b^2=0$, $ab=ba$ and $bc=cb$. The black points indicate a tame concealed algebra of type $\widetilde{\mathbb{E}}_6$ (see No.1 in \cite{B-critical}). Hence, $N(3)\otimes\mathbf{B}_4$ is representation-infinite.
\end{proof}

Proposition \ref{prop::N3-B4} implies that if $N(n)\otimes B$ is representation-finite, then the quiver $Q_B$ of $B$ must be of type $A$ or $\widetilde{A}$. Note that $B$ is simply connected. If $|B|\ge 4$ and $Q_B$ is of type $\widetilde{A}$, then $B$ is isomorphic to $\Omega_{s,t}$ with $s+t=|B|$, where $\Omega_{s,t}$ is defined by
$$
\vcenter{\xymatrix@C=1cm@R=0.1cm{
&2 \ar[r]&3 \ar[r]&\cdots \ar[r] &s\ar[dr] &\\
1 \ar[ur] \ar[dr]\ar@{.}[rrrrr] &&&&&s+t \\
&s+1 \ar[r]  &s+2 \ar[r] &\cdots \ar[r] &s+t-1\ar[ur] &}}
$$
with the unique commutativity relation.

\begin{proposition}\label{prop::N(n)-A_st}
If $s+t\ge 4$, then $N(n)\otimes \Omega_{s,t}$ is representation-infinite for any $n\geq 3$.
\end{proposition}
\begin{proof}
It follows from Proposition \ref{prop::A2} that $\mathbf{A}_2^{(+)}\otimes \Omega_{s,t}$ is representation-infinite. Since there is a surjection from $N(n)$ to $\mathbf{A}_2^{(+)}$, we obtained the statement.
\end{proof}

By combining Proposition \ref{prop::N3-B4} and Proposition \ref{prop::N(n)-A_st}, the quiver $Q_B$ of $B$ should be of type $A$ if $N(n)\otimes B$ is representation-finite. This reduces the problem to some small cases.

\begin{proposition}\label{prop::N(n)-B5}
Suppose $|B|=4$ and $Q_B$ is of the form $\vcenter{\xymatrix@C=0.5cm@R=0.5cm{
1\ar@{-}[r]&2\ar@{-}[r]&3\ar@{-}[r]&4}}$. Then, $N(n)\otimes B$ with $n\ge 3$ is representation-finite if and only if $B$ or $B^{\mathsf{op}}$ is isomorphic to
\begin{equation}
\mathbf{B}_5:=\mathsf{k}\left(
\xymatrix@C=1cm@R=0.5cm{
1\ar[r]^\alpha&2&3\ar[l]_\beta&4\ar[l]_\gamma}
\right)
\bigg/
\langle\gamma\beta\rangle.
\end{equation}
\end{proposition}
\begin{proof}
Since $B$ is not Nakayama, $B$ cannot isomorphic to $\mathbf{A}_4^{(+++)}$ or $\mathbf{A}_4^{(---)}$.
If $N(n)\otimes B$ is representation-finite, then $B \not \simeq \mathbf{A}_4^{\varepsilon}$ for any choice of $\varepsilon$ (see Lemma \ref{lem::A4}) and hence, the unique possible case is $B\simeq \mathbf{B}_5$ (or $\mathbf{B}_5^{\mathsf{op}}$). In fact, $N(n)\otimes \mathbf{B}_5$ is representation-finite following from \cite[Page 155, (A)]{L-rep-type}.
\end{proof}

We define
\begin{equation}
\mathbf{B}_6:=\mathsf{k}\left(
\xymatrix@C=1cm@R=0.5cm{
1\ar[r]^\alpha&2&3\ar[l]_\beta&4\ar[l]_\gamma\ar[r]^\delta&5}
\right)
\bigg/
\langle\gamma\beta\rangle,
\end{equation}
and
\begin{equation}
\mathbf{B}_7:=\mathsf{k}\left(
\xymatrix@C=1cm@R=0.5cm{
5\ar[r]^\delta &1\ar[r]^\alpha&2&3\ar[l]_\beta&4\ar[l]_\gamma}
\right)
\bigg/
\langle\gamma\beta, \delta\alpha\rangle.
\end{equation}

\begin{lemma}\label{lem::N3-B6-B7}
If $B\simeq \mathbf{B}_6$ or $\mathbf{B}_7$, then $N(3)\otimes B$ is representation-infinite.
\end{lemma}
\begin{proof}
The tensor product $N(3)\otimes \mathbf{B}_6$ is given by the quiver
$$
\vcenter{\xymatrix@C=1.5cm@R=0.7cm{
\circ \ar[r]^a\ar[d]^b\ar@{.}[dr]&\bullet \ar[d]^b &\circ\ar[l]_c\ar[d]^b\ar@{.}[dl]&\circ\ar[l]_c \ar[r]^d\ar[d]^b\ar@{.}[dr]\ar@{.}[dl]&\circ\ar[d]^b\\
\bullet\ar[r]^a\ar[d]^b\ar@{.}[dr]&\bullet\ar[d]^b &\bullet\ar[l]_c\ar[d]^b\ar@{.}[dl]&\circ\ar[l]_c \ar[r]^d\ar[d]^b\ar@{.}[dr]\ar@{.}[dl]&\bullet\ar[d]^b\\
\bullet \ar[r]^a&\circ  &\bullet\ar[l]_c&\bullet\ar[l]_c \ar[r]^d&\bullet
}}
$$
with $b^2=c^2=0$, $ab=ba$, $cb=bc$ and $db=bd$. In the quiver, the black points induce a tame concealed algebra of type $\widetilde{\mathbb{E}}_8$ (see No.4 in \cite{B-critical}) and thus, $N(3)\otimes \mathbf{B}_6$ is representation-infinite. Similarly, $N(3)\otimes \mathbf{B}_7$ is given by the quiver
$$
\vcenter{\xymatrix@C=1.5cm@R=0.7cm{
\circ \ar[r]^a\ar[d]^b\ar@{.}[dr]&\circ\ar@{.}[dr]\ar[r]^a\ar[d]^b &\bullet\ar[d]^b&\circ\ar[l]_c \ar[d]^b\ar@{.}[dl]&\circ\ar[l]_c\ar[d]^b\ar@{.}[dl]\\
\circ \ar[r]^a\ar[d]^b\ar@{.}[dr]&\bullet \ar@{.}[dr]\ar[r]^a\ar[d]^b &\bullet\ar[d]^b&\bullet\ar[l]_c \ar[d]^b\ar@{.}[dl]&\circ\ar[l]_c\ar[d]^b\ar@{.}[dl]\\
\bullet\ar[r]^a&\bullet \ar[r]^a&\circ&\bullet\ar[l]_c &\bullet\ar[l]_c
}}
$$
with $a^2=b^2=c^2=0$, $ab=ba$ and $cb=bc$. Here, the black points give a tame concealed algebra of type $\widetilde{\mathbb{E}}_7$ (see No.2 in \cite{B-critical}).
\end{proof}

\begin{proposition}\label{prop::N(3)-B-not-naka}
Suppose $|B|\ge 5$ and $Q_B$ is of type $A$. Then, the tensor product $N(3)\otimes B$ is representation-finite if and only if $B$ does not contain one of $\mathbf{A}_4^{\varepsilon}$ for some $\varepsilon$, $\mathbf{B}_6$, $\mathbf{B}_7$ and their opposite algebras as a quotient algebra.
\end{proposition}
\begin{proof}
We assume that $B$ does not contain one of $\mathbf{A}_4^{\varepsilon}$, $\mathbf{B}_6$, $\mathbf{B}_7$ and their opposite algebras as a quotient. Since $B$ is not a Nakayama algebra, the quiver $Q_B$ of $B$ contains at least one $\xymatrix@C=0.6cm{\circ&\circ\ar[l]\ar[r]&\circ}$ or $\xymatrix@C=0.6cm{\circ\ar[r]&\circ &\circ\ar[l]}$ as a subquiver. If $Q_B$ contains
$$
\xymatrix@C=0.6cm{\circ\ar@{-}[r]&\circ&\circ\ar[l]\ar[r]&\circ\ar@{-}[r]&\circ}  \quad\text{or}\quad \xymatrix@C=0.6cm{\circ\ar@{-}[r]&\circ\ar[r]&\circ &\circ\ar[l]\ar@{-}[r]&\circ}
$$
as a subquiver, this contradicts our assumption. Thus, $\xymatrix@C=0.6cm{\circ&\circ\ar[l]\ar[r]&\circ}$ and $\xymatrix@C=0.6cm{\circ\ar[r]&\circ &\circ\ar[l]}$ can only appear in the leftmost or rightmost of $Q_B$. Moreover, $Q_B$ cannot be of the form
$$
\xymatrix@C=0.6cm{\circ\ar[r]&\circ&\circ\ar[l]\ar@{-}[r]&\circ\ar@{-}[r]&\cdots\ar@{-}[r]&\circ\ar@{-}[r]&\circ\ar[r]&\circ&\circ\ar[l]}
$$
or
$$
\xymatrix@C=0.6cm{\circ&\circ\ar[r]\ar[l]&\circ\ar@{-}[r]&\circ\ar@{-}[r]&\cdots\ar@{-}[r]&\circ\ar@{-}[r]&\circ&\circ\ar[r]\ar[l]&\circ},
$$
otherwise, we will also obtain a contradiction. In summary, $Q_B$ should be of the form
$$
\xymatrix@C=0.6cm{\circ\ar[r]&\circ&\circ\ar[l]&\circ\ar[l]&\cdots\ar[l]&\circ\ar[l]&\circ\ar[l]&\circ\ar[l]&\circ\ar[l]}
$$
or
$$
\xymatrix@C=0.6cm{\circ\ar[r]&\circ&\circ\ar[l]&\circ\ar[l]&\cdots\ar[l]&\circ\ar[l]&\circ\ar[l]&\circ\ar[r]\ar[l]&\circ},
$$
up to isomorphism.

Similar to the proof of Proposition \ref{prop::N3-Nakayama}, we show that each $3\times 8$ truncation (if it exists) of $N(3)\otimes B$ is representation-finite. Up to isomorphism, a $3\times 8$ truncation of $N(3)\otimes B$ is isomorphic to one of the algebras listed below:
\begin{enumerate}
\item $\xymatrix@C=0.6cm{1\ar[r]^{a_1}&2&3\ar[l]_{a_2}&4\ar[l]_{a_3}&5\ar[l]_{a_4}&6\ar[l]_{a_5}&7\ar[l]_{a_6}&8\ar[l]_{a_7}}$ with $a_3a_2=a_5a_4=a_7a_6=0$.

\item [($1^1$)] $\xymatrix@C=0.6cm{1\ar[r]^{a_1}&2&3\ar[l]_{a_2}&4\ar[l]_{a_3}&5\ar[l]_{a_4}&6\ar[l]_{a_5}&7\ar[l]_{a_6}&8\ar[l]_{a_7}}$ with $a_3a_2=a_4a_3=a_5a_4=a_7a_6=0$.

\item $\xymatrix@C=0.6cm{1\ar[r]^{a_1}&2&3\ar[l]_{a_2}&4\ar[l]_{a_3}&5\ar[l]_{a_4}&6\ar[l]_{a_5}&7\ar[l]_{a_6}&8\ar[l]_{a_7}}$ with $a_3a_2=a_5a_4=a_6a_5=0$.

\item[($2^1$)] $\xymatrix@C=0.6cm{1\ar[r]^{a_1}&2&3\ar[l]_{a_2}&4\ar[l]_{a_3}&5\ar[l]_{a_4}&6\ar[l]_{a_5}&7\ar[l]_{a_6}&8\ar[l]_{a_7}}$ with $a_3a_2=a_5a_4=a_6a_5=a_7a_6=0$.

\item $\xymatrix@C=0.6cm{1\ar[r]^{a_1}&2&3\ar[l]_{a_2}&4\ar[l]_{a_3}&5\ar[l]_{a_4}&6\ar[l]_{a_5}&7\ar[l]_{a_6}&8\ar[l]_{a_7}}$ with $a_3a_2=a_4a_3=a_6a_5=0$.

\item [($3^1$)] $\xymatrix@C=0.6cm{1\ar[r]^{a_1}&2&3\ar[l]_{a_2}&4\ar[l]_{a_3}&5\ar[l]_{a_4}&6\ar[l]_{a_5}&7\ar[l]_{a_6}&8\ar[l]_{a_7}}$ with $a_3a_2=a_4a_3=a_6a_5=a_7a_6=0$.

\item [($3^2$)] $\xymatrix@C=0.6cm{1\ar[r]^{a_1}&2&3\ar[l]_{a_2}&4\ar[l]_{a_3}&5\ar[l]_{a_4}&6\ar[l]_{a_5}&7\ar[l]_{a_6}&8\ar[l]_{a_7}}$ with $a_3a_2=a_4a_3=a_5a_4=a_6a_5=0$.

\item [($3^3$)] $\xymatrix@C=0.6cm{1\ar[r]^{a_1}&2&3\ar[l]_{a_2}&4\ar[l]_{a_3}&5\ar[l]_{a_4}&6\ar[l]_{a_5}&7\ar[l]_{a_6}&8\ar[l]_{a_7}}$ with $a_{i+1}a_i=0$ for $i=2,3,4,5,6$.

\item $\xymatrix@C=0.6cm{1\ar[r]^{a_1}&2&3\ar[l]_{a_2}&4\ar[l]_{a_3}&5\ar[l]_{a_4}&6\ar[l]_{a_5}&7\ar[l]_{a_6}\ar[r]^{a_7}&8}$ with $a_3a_2=a_5a_4=a_6a_5=0$.

\item $\xymatrix@C=0.6cm{1\ar[r]^{a_1}&2&3\ar[l]_{a_2}&4\ar[l]_{a_3}&5\ar[l]_{a_4}&6\ar[l]_{a_5}&7\ar[l]_{a_6}\ar[r]^{a_7}&8}$ with $a_3a_2=a_4a_3=a_6a_5=0$.

\item [($5^1$)] $\xymatrix@C=0.6cm{1\ar[r]^{a_1}&2&3\ar[l]_{a_2}&4\ar[l]_{a_3}&5\ar[l]_{a_4}&6\ar[l]_{a_5}&7\ar[l]_{a_6}\ar[r]^{a_7}&8}$ with $a_3a_2=a_4a_3=a_5a_4=a_6a_5=0$.
\end{enumerate}
In the list, the algebras $(i^j)$'s are quotient algebras of ($i$). Then, one may use the software GAP 4.10.2 to check that the Tits form of $N(3)\otimes B$ in each case is weakly positive, we omit the details. It turns out that $N(3)\otimes B$ is representation-finite.

If $N(3)\otimes B$ does not have a $3\times 8$ truncation, it suffices to consider the following cases:
\begin{enumerate}
\item $\xymatrix@C=0.6cm{1\ar[r]^{a_1}&2&3\ar[l]_{a_2}&4\ar[l]_{a_3}&5\ar[l]_{a_4}\ar[r]^{a_5}&6}$ with $a_3a_2=a_4a_3=0$.

\item $\xymatrix@C=0.6cm{1\ar[r]^{a_1}&2&3\ar[l]_{a_2}&4\ar[l]_{a_3}&5\ar[l]_{a_4}&6\ar[l]_{a_5}\ar[r]^{a_6}&7}$ with $a_3a_2=a_5a_4=0$.

\item [($2^1$)] $\xymatrix@C=0.6cm{1\ar[r]^{a_1}&2&3\ar[l]_{a_2}&4\ar[l]_{a_3}&5\ar[l]_{a_4}&6\ar[l]_{a_5}\ar[r]^{a_6}&7}$ with $a_3a_2=a_4a_3=a_5a_4=0$.
\end{enumerate}
It is not difficult to see that $N(3)\otimes B$ is representation-finite in all three cases above.
\end{proof}

\begin{proposition}\label{prop::N4-B8}
Suppose $|B|\ge 5$ and $Q_B$ is of type $A$. Then, the tensor product $N(n)\otimes B$ is representation-infinite for any $n\geq 4$.
\end{proposition}
\begin{proof}
By Proposition \ref{prop::N(3)-B-not-naka}, it suffices to show that $N(4)\otimes \mathbf{B}_8$ is representation-infinite, where
\begin{equation}
\mathbf{B}_8:=\mathsf{k}\left(
\xymatrix@C=1cm@R=0.5cm{
1\ar[r]^\alpha&2&3\ar[l]_\beta&4\ar[l]_\gamma&5\ar[l]_\delta}
\right)
\bigg/\langle\gamma\beta, \delta\gamma \rangle.
\end{equation}
In fact, $N(4)\otimes \mathbf{B}_8$ is presented by the quiver
$$
\vcenter{\xymatrix@C=1.5cm@R=0.7cm{
\circ \ar[r]^a\ar[d]^b\ar@{.}[dr]&\bullet\ar[d]^b &\circ\ar[d]^b\ar[l]_c\ar@{.}[dl]&\circ\ar[l]_c \ar[d]^b\ar@{.}[dl]&\circ\ar[l]_c\ar[d]^b\ar@{.}[dl]\\
\bullet \ar[r]^a\ar[d]^b\ar@{.}[dr]&\bullet \ar[d]^b &\bullet\ar@{.}[dl]\ar[d]^b\ar[l]_c&\circ\ar[l]_c \ar[d]^b\ar@{.}[dl]&\circ\ar[l]_c\ar[d]^b\ar@{.}[dl]\\
\bullet \ar[r]^a\ar[d]^b\ar@{.}[dr]&\circ \ar[d]^b &\bullet\ar@{.}[dl]\ar[d]^b\ar[l]_c&\bullet\ar[l]_c \ar[d]^b\ar@{.}[dl]&\circ\ar[l]_c\ar[d]^b\ar@{.}[dl]\\
\circ\ar[r]^a&\circ&\circ\ar[l]_c&\bullet\ar[l]_c &\bullet\ar[l]_c
}}
$$
with $b^2=c^2=0$, $ab=ba$ and $bc=cb$. It is obvious that the black points in the quiver give a tame concealed algebra of type $\widetilde{\mathbb{E}}_8$ (see No.4 in \cite{B-critical}).
\end{proof}

In summary, we have the following result related to our classification.
\begin{theorem}\label{theo::naka-rad-zeor-not-naka}
Let $A\simeq N(n)$ with $n\ge 3$ and $B$ a simply connected but not Nakayama algebra with $|B|\ge 3$. Then, $A\otimes B$ is representation-finite if and only if one of the following holds:
\begin{enumerate}
\item[\rm{(1)}] $|B|=3$.
\item[\rm{(2)}] $|B|=4$ and $B$ is isomorphic to $\mathbf{B}_5$ or $\mathbf{B}_5^{\textsf{op}}$.
\item[\rm{(3)}] $|B|\ge 5$, $n=3$ and $B$ does not contain one of $\mathbf{A}_4^{\varepsilon}$ for some $\varepsilon$, $\mathbf{B}_6$, $\mathbf{B}_7$ and their opposite algebras as a quotient algebra.
\end{enumerate}
\end{theorem}

\section{Minimal representation-infinite cases}
In Section 3 and Section 4, we have obtained a complete list of \emph{minimal representation-infinite} simply connected tensor products, in the sense that $A\otimes B$ is representation-infinite, but any proper quotient $A'\otimes B'$ is representation-finite, for $A'$ and $B'$ being quotients of $A$ and $B$ respectively. More precisely, the list is given by 
\begin{enumerate}
\item $A\simeq \mathbf{A}_2^{(+)}$ and $B$ (or $B^\mathsf{op}$) is isomorphic to one of the algebras in the list (IT) in \cite{LS-tame-triangular-matrix}, see Proposition \ref{prop::A2}.

\item $A\simeq \mathbf{A}_3^{\varepsilon}$ and $B\simeq \mathbf{A}_3^{\omega}$ for some $\varepsilon$ and $\omega$, see Lemma \ref{lem::A3-A3}.

\item $A\simeq \mathbf{A}_4^{\varepsilon}$ for some $\varepsilon$, and $B\simeq N(3)$, see Lemma \ref{lem::A4} and Lemma \ref{lem::N3-A4}.

\item $A\simeq N(3)$ and $B$ is isomorphic to one of $\mathbf{B}_i$, $\mathbf{B}_i^{\mathsf{op}}$ for $i\in \{1,2,4,6,7\}$, see Lemma \ref{lem::N3-B1}, Lemma \ref{lem::N3-B2}, Proposition \ref{prop::N3-B4} and Lemma \ref{lem::N3-B6-B7}.

\item $A\simeq N(4)$ and $B$ is isomorphic to one of $\mathbf{B}_3$, $\mathbf{B}_8$, $\mathbf{B}_3^{\mathsf{op}}$, $\mathbf{B}_8^{\mathsf{op}}$, see Lemma \ref{lem::N4-B3} and Proposition \ref{prop::N4-B8}.
\end{enumerate}

\ \\
\section*{Acknowledgements}
The authors are grateful to Takahide Adachi, Takuma Aihara and Takahiro Honma for many useful conversations and suggestions.
KM is partially supported by JSPS Grant-in-Aid for Young Scientists (Grant No. 20K14302). QW is partially supported by JSPS Grant-in-Aid for Young Scientists (Grant No. 20J10492), National Key Research and Development Program of China (Grant No. 2020YFA0713000) and China Postdoctoral Science Foundation (Grant No. YJ20220119 and No. 2023M731988).

\ \\

\begin{thebibliography}{AAAA}
\setlength{\baselineskip}{14.2pt}
\bibitem[AIR]{AIR}
{T. Adachi, O. Iyama and I. Reiten},
$\tau$-tilting theory.
{\it Compos. Math.} {\bf 150} (2014), no. 3, 415--452.

\bibitem[AH]{Aihara-Honma}
T. Aihara and T. Honma,
$\tau$-tilting finite triangular matrix algebras.
{\it J. Pure Appl. algebra} {\bf 225} (2021), no 12, 106785.

\bibitem[A]{Assem-simply-connected}
I. Assem,
Simply connected algebras.
{\it Resenhas IME-USP} {\bf 4} (1999), no. 2, 93--125.

\bibitem[ABM]{ABM}
I. Assem, J. C. Bustamante and P. L. Meur,
Special biserial algebras with no outer derivations.
{\it Colloq. Math.} {\bf 125} (2011), no. 1, 83--98.

\bibitem[AR]{AR-triangular-matrix}
M. Auslander and I. Reiten,
On the representation type of triangular matrix rings.
{\it J. London Math. Soc. (2)} {\bf 12} (1976), 371--382.

\bibitem[ARS]{ARS}
M. Auslander, I. Reiten and O. Smal\o,
Representation Theory of Artin algebras.
{\it Cambridge Studies in Advanced Mathematics} \textbf{36}, {\it Cambridge University Press}, 1995.

\bibitem[AS]{AS}
I. Assem and A. Skowro$\acute{\text{n}}$ski,
On some classes of simply connected algebras.
{\it Proc. London Math. Soc.} {\bf 56} (1988), no. 3, 417--450.\

\bibitem[ASS]{ASS}
I. Assem, D. Simson and A. Skowro$\acute{\text{n}}$ski,
Elements of the Representation Theory of Associative Algebras, vol 1, Techniques of Representation Theory,
{\it London Mathematical Society Student Texts} \textbf{65}, {\it Cambridge University Press}, 2006.

\bibitem[Bo1]{B-tits-form}
K. Bongartz,
Algebras and quadratic forms.
{\it J. London Math. Soc. (2)} {\bf 28} (1983), no. 3, 461--469.

\bibitem[Bo2]{B-critical}
K. Bongartz,
Critical simply connected algebras.
{\it Manuscripta Math.} {\bf 46} (1984), no. 1-3, 117--136.

\bibitem[BD]{BD-finite-group}
V. M. Bondar enko and Yu. A. Drozd,
The representation type of finite groups, in: Modules and Representations.
{\it Zap. Nauchn. Sem. LOMI} {\bf 57} (1977), 24--41. (in Russian)

\bibitem[BoG]{BG-covering}
K. Bongartz and P. Gabriel,
Covering spaces in representation theory.
{\it Invent. Math.} {\bf 65} (1982), no. 3, 331--378.

\bibitem[BrG]{BG-standard form}
O. Bretscher and P. Gabriel,
The standard form of a representation-finite algebra.
{\it Bull. Soc. Math. France} {\bf 111} (1983), no. 1, 21--40.

\bibitem[D]{Dr-tame-wild}
Yu. A. Drozd,
Tame and wild matrix problems, in: Representation Theory II, Lecture Notes in Math., vol. 832.
{\it Springer Verlag}, (1980), pp. 242--258.

\bibitem[HM]{HM-selfinjective}
M. Hoshino and I. Miyachi,
Tame triangular matrix algebras over self-injective algebras.
{\it Tsukuba J. Math.} {\bf 11} (1987), 383--391.

\bibitem[HR]{HR-tilted}
D. Happel and C.M. Ringel,
Tilted algebras.
{\it Trans. Amer. Math. Soc.} {\bf 274} (1982), no. 2, 399--443.

\bibitem[HV]{HV}
D. Happel and D. Vossieck,
Minimal algebras of infinite representation type with preprojective component.
{\it Manuscripta Math.} {\bf 42} (1983), 221--243.

\bibitem[L1]{L-rep-type}
Z. Leszczy$\acute{\text{n}}$ski,
On the representation type of tensor product algebras.
{\it Fund. Math.} {\bf 144} (1994), no. 2, 143--161.

\bibitem[L2]{L-special-alg}
Z. Leszczy$\acute{\text{n}}$ski,
On the representation type of triangular matrix algebras over special algebras.
{\it Fund. Math.} {\bf 137} (1991), 65--80.

\bibitem[LS1]{LS-tame-triangular-matrix}
Z. Leszczy$\acute{\text{n}}$ski and A. Skowro$\acute{\text{n}}$ski,
Tame triangular matrix algebras.
{\it Colloq. Math.} {\bf 86} (2000), no. 2, 259--303.

\bibitem[LS2]{LS-tame-tensor-product}
Z. Leszczy$\acute{\text{n}}$ski and A. Skowro$\acute{\text{n}}$ski,
Tame tensor product algebras.
{\it Colloq. Math.} {\bf 98} (2003), no. 1, 125--145.

\bibitem[MP]{MP}
R. Mart\'{i}nez-Villa and J. A. de la Pe$\tilde{\text{n}}$a,
The universal cover of a quiver with relations.
{\it J. Pure Appl. algebra}  {\bf 30}  (1983), 277--292.

\bibitem[SS]{SS2}
D. Simson and A. Skowro$\acute{\text{n}}$ski,
Elements of the Representation Theory of Associative Algebras, vol 3, Representation-Infinite Tilted Algebras,
{\it London Mathematical Society Student Texts} {\bf 72}, {\it Cambridge University Press}, 2007.

\bibitem[S]{S-nakayama}
A. Skowro$\acute{\text{n}}$ski,
Tame triangular matrix algebras over Nakayama algebras.
{\it J. London Math. Soc.} {\bf 34} (1986), 245--264.

\bibitem[W]{W-simply}
Q. Wang,
On $\tau$-tilting finite simply connected algebras.
{\it Tsukuba J. Math.} {\bf (1) 46} (2022), 1--37.
\end{thebibliography}
\end{document}